\documentclass{article}

\usepackage{amsfonts}
\usepackage{amsmath}
\usepackage{amssymb}
\usepackage{theorem}
\usepackage{latexsym}

\numberwithin{equation}{section}

\newtheorem{thm}{Theorem}[section]
\newtheorem{co}[thm]{Corollary}
\newtheorem{prop}[thm]{Proposition}
\newtheorem{lemma}[thm]{Lemma}
\newtheorem{rmk}[thm]{Remark}
\newtheorem{eg}[thm]{Example}
\newtheorem{Defn}[thm]{Definition}

\newcommand{\eosp}{\ensuremath{\operatorname{eosp}}}
\newcommand{\osp}{\ensuremath{\operatorname{osp}}}
\newcommand{\ospfd}{\ensuremath{\operatorname{osp_{\,fd}}}}
\newcommand{\ca}{\ensuremath{{\cal{A}}}}
\newcommand{\cas}{\ca\ }
\newcommand{\cu}{\ensuremath{{\cal{U}}}}
\newcommand{\cus}{\cu\ }
\newcommand{\es}{\ensuremath{{\cal{S}}}}
\newcommand{\ess}{\es\ }
\newcommand{\esone}{\ensuremath{{\cal{S}}^{(1)}}}

\newcommand{\estwo}{\ensuremath{{\cal{S}}^{(2)}}}

\newcommand{\esthree}{\ensuremath{{\cal{S}}^{(3)}}}

\newcommand{\te}{\ensuremath{{\cal{T}}}}
\newcommand{\tes}{\te\ }
\newcommand{\teone}{\ensuremath{{\cal{T}}^{(1)}}}

\newcommand{\tetwo}{\ensuremath{{\cal{T}}^{(2)}}}

\newcommand{\tethree}{\ensuremath{{\cal{T}}^{(3)}}}

\newcommand{\m}{\ensuremath{{\cal{M}}}}
\newcommand{\ms}{\m\ }
\newcommand{\n}{\ensuremath{{\cal{N}}}}
\newcommand{\ns}{\n\ }
\newcommand{\p}{\ensuremath{{\cal{P}}}}
\newcommand{\ps}{\p\ }
\newcommand{\calr}{\ensuremath{{\cal{R}}}}
\newcommand{\calrs}{\calr\ }
\newcommand{\de}{\ensuremath{{\cal{D}}}}

\newcommand{\ce}{\ensuremath{{\cal{E}}}}
\newcommand{\ces}{\ce\ }
\newcommand{\g}{\ensuremath{{\cal{G}}}}

\newcommand{\und}[1]{\underline{{#1\!}}\,}
\newcommand{\ovl}[1]{\overline{{#1}}}
\newcommand{\qed}{\hfill\rule{2mm}{2mm}}

\newcommand{\ve}{\ensuremath{\, |\, }}
\newcommand{\half}{\ensuremath{\frac{1}{2}}}
\newcommand{\til}{\ensuremath{\tilde{\ }}}
\newcommand{\sone}{\textup{(S1)}}
\newcommand{\sones}{\sone\ }
\newcommand{\stwo}{\textup{(S2)}}
\newcommand{\stwos}{\stwo\ }
\newcommand{\sthree}{\textup{(S3)}}
\newcommand{\sthrees}{\text{\sthree}\ }
\newcommand{\tone}{\textup{(T1)}}

\newcommand{\tones}{\tone\ }
\newcommand{\ttwo}{\textup{(T2)}}
\newcommand{\ttwos}{\text{\ttwo}\ }
\newcommand{\tthree}{\textup{(T3)}}

\newcommand{\inv}{\ensuremath{^{-1}}}
\newcommand{\idd}{\ensuremath{Id}}
\newcommand{\C}{\ensuremath{\mathbb{C}}}
\newcommand{\Z}{\ensuremath{\mathbb{Z}}}
\newcommand{\N}{\ensuremath{\mathbb{N}}}
\newcommand{\F}{\ensuremath{\mathbb{F}}}
\newcommand{\e}{\ensuremath{\mathbf{E}}}
\newcommand{\id}{\ensuremath{\!\cdot\! \idd}}
\newcommand{\ad}{\operatorname{\bf{ad}}}
\newcommand{\adloc}{\ensuremath{\operatorname{\bf{ad_{loc}}}}}
\newcommand{\ba}{\ensuremath{A}}
\newcommand{\bas}{\ba\ }
\newcommand{\lk}{\ensuremath{k}}

\newcommand{\bk}{\ensuremath{K}}
\newcommand{\bks}{\bk\ }
\newcommand{\bastar}{\ensuremath{A^{\times}}}
\newcommand{\bastars}{\bastar\ }
\newcommand{\bz}{\ensuremath{\operatorname{\bar{0}}}}
\newcommand{\bo}{\ensuremath{\operatorname{\bar{1}}}}
\newcommand{\q}{\ensuremath{q}}
\newcommand{\qs}{\q\ }
\newcommand{\der}{\ensuremath\operatorname{Der}}
\newcommand{\ider}{\ensuremath{\ad}}
\newcommand{\ed}{\ensuremath\operatorname{End}}
\newcommand{\hm}{\ensuremath\operatorname{Hom}}
\newcommand{\ann}{\ensuremath\operatorname{Ann}}

\begin{document}
\rmfamily
\bibliographystyle{unsrt}

\title{{DERIVATIONS OF ORTHOSYMPECTIC LIE
SUPERALGEBRAS}}
\author{Ana Duff\footnote{The author gratefully acknowledges the support
of National Science and Engineering Research Council Postgraduate
Scholarships PGS A and PGS B and the support of the Ontario
Graduate Scholarship in Science and Technology.}\ \footnote{This
paper is part of the author's Ph.D. dissertation \cite{D} written
under the supervision of E. Neher at the University of Ottawa. The
author wishes to thank him for his support and
guidance.}\\
\ \\
Department of Mathematics and Statistics, York University,\\
Toronto, ON M3J 1P3, Canada\\aduff@mathstat.yorku.ca}
\date{}\maketitle
\begin{center}{\bf ABSTRACT}\\
\ \\
\begin{quote}In this paper we describe the derivations of orthosymplectic
Lie superalgebras over a superring.  In particular, we derive
sufficient conditions under which the derivations can be expressed
as a semidirect product of inner and outer derivations.  We then
present some examples for which these conditions hold.\end{quote}
\end{center}
\ \\

\begin{center}{\bf  INTRODUCTION}\end{center}
\ \\
The original motivation for this paper came from the study of
extended affine Lie algebras (\cite{AABGP}).  In particular we
were interested in describing the derivations of Lie algebras over
fields of characteristic 0 graded by root systems of type $B_I$
and $D_I$ and developing conditions under which we could write the
derivation algebra as a semidirect product of the inner and the
outer derivations.  Such a decomposition is very useful in the
construction of extended affine Lie algebras of type $B$ and $D$
from their centreless core (\cite{AG}).  In light of the recent
activity in developing a theory of root-graded Lie superalgebras
(\cite{BeE}, \cite{GN}), it is of interest to consider Lie
superalgebras.  Our methods will allow us to describe derivations
of various types of orthosymplectic Lie superalgebra over
superrings containing \half,
including some examples which occur in the recent paper \cite{CW}.\\
\\
Some work on derivations of Lie superalgebras was done by Kac in
\cite{Kac2} and Scheunert in \cite{Scheunert}.  In particular,
they described the derivations of simple finite-dimensional Lie
superalgebras over algebraically closed fields of characteristic
0.  For orthosymplectic Lie superalgebras, their result is a
special case of Corollary \ref{kac and scheunert} presented in
this paper.  Significant work on derivations was recently done by
Benkart in \cite{Be}.  She described the derivations of Lie
algebras over fields of characteristic zero which are graded by
finite root systems using the derivations of the coordinate
algebras.  The result described in Corollary \ref{S for almost
diagonalizable+P} is a generalization of Benkart's result when we
consider the Lie algebras which are graded by root systems of type
$B_I,\,|I|\geq 3$ and $D_I,\,|I|\geq 4$.

We assume that \bks is a supercommutative and unital superring
containing \half\ and \bas a superextension of \bk.  We consider
subalgebras of the orthosymplectic Lie superalgebra \osp(\q) where
\qs is an \ba-quadratic form on an \ba-supermodule \m.  The
superalgebra \osp(\q) is the Lie superalgebra of all
\ba-endomorphisms $x$ of \ms such that $q(x(m),n)+(-
1)^{|m||n|}q(x(n),m)=0$ for all $m,n\in \m$.  It has an ideal,
\eosp(\q), which is defined to be the \Z-span of the maps
$\e_{m,n}$ for homogeneous $m,n\in\m$ where
$\e_{m,n}(p)=mq(n,p)-(- 1)^{|n||p|}q(m,p)n$ for $p\in\m$.  We show
that if $\q_{\infty}$ is the orthogonal sum of the hyperbolic
superplane and another quadratic form \qs on an \ba-supermodule
\m, then the \bk-derivations of any subalgebra $\ce_{\infty}$ of
$\osp(\q_{\infty})$ containing $\eosp(\q_{\infty})$ can be
described as a sum of the inner derivations and a Lie superalgebra
$\es\oplus\te$ for certain $\es,\,\te\subset\ed_K\m$.  We also
determine the intersection of the inner derivations and the
superalgebra $\es\oplus\te$ and we determine the conditions under
which we can write the algebra of derivations as a semidirect
product of the inner derivations and a certain subalgebra.
Finally, we describe some examples where we do get the splitting
of \bk-derivations of $\ce_{\infty}$, $\der_{\bk}(\ce_{\infty})$,
into a semidirect product of the inner derivations and a
subalgebra of $\es\oplus\te$.  In particular, when \qs is an
almost diagonalizable \ba-quadratic form on a free \ba-supermodule
of dimension greater than $ 2$ (i.e., there exists a homogeneous
basis $\{m_i\mid i\in I\}$ of \ms such that for each $i\in I$
there exists $\und{i}\in I$ such that $q(m_i,m_j)$ is 0 if $i\neq
\und{i}$ and is invertible otherwise), we get the semidirect
splitting of the derivations for a number of subalgebras
$\ce_{\infty}$.  These include the centreless core $L$ of an
extended affine Lie algebra over \C\ of type $B_I$ or $D_I$ and
nullity $\nu$ and, in this case in particular, we get
$\der_{\C}(L/Z(L))\cong\ider (L/Z(L))\rtimes\der_{\C}\C[t_1^{\pm
1},\ldots,t_{\nu}^{\pm 1}]$ for $I$ a finite set.


\section{SUPERALGEBRAS}
In this section we will describe some fundamental concepts of
su\-per\-struc\-tur\-es.  The reader can find more extensive
coverage of this material in \cite{Kac2}, \cite{Scheunert} and
\cite{GN}.\newline\newline Let $\Z_2=\{\bz,\bo\}$ be the field of
two elements.  A ring $K$ is called a \emph{superring} if $K$ is a
$\Z_2$-graded ring, i.e., $K=K_{\bz}\oplus K_{\bo}$ such that
$K_{\alpha}K_{\beta}\in K_{\alpha+\beta}$ for all
$\alpha,\beta\in\Z_2$. We say that a superring \bks is
\emph{supercommutative} if $ab=(-1)^{\alpha\beta}ba$ for all
$a\in\bk_{\alpha},b\in K_{\beta}$ for all $\alpha,\beta\in\Z_2$.
Throughout this paper \bks shall denote a unital and
supercommutative superring and we will assume that \bks contains
\half.  Note that $1\in K_{\bz}$.
\newline\newline
A \bk-bimodule \ms is called a \emph{\bk-supermodule} if
$\m=\m_{\bz}\oplus\m_{\bo}$ for some $\bk_{\bz}$\,-\,submodules
$\m_{\bz}$ and $\m_{\bo}$ of \ms such that
$K_{\alpha}\m_{\beta}\subset\m_{\alpha+\beta}$ for all
$\alpha,\beta\in\Z_2$ and the \bk-module action satisfies
\begin{eqnarray}
\label{supermodule} am=(-1)^{\alpha\beta}ma
\end{eqnarray}
for all $a\in\bk_{\alpha},\,m\in\m_{\beta}$.  Throughout this
paper, any time we talk about a $\Z_2$-graded structure
$Z=Z_{\bz}\oplus Z_{\bo}$ we will call an element $z\in Z$
\emph{homogeneous} if $z\in Z_{\alpha}$ for some $\alpha\in\Z_2$
and we will say that $z$ is of \emph{degree} $\alpha$, denoted by
$|z|=\alpha$; in this case we will assume throughout that when
$|z|$ occurs in an expression, then it is assumed that $z$ is
homogeneous, and that the expression extends to the other elements
by linearity.  When we refer to a submodule of a supermodule we
assume that the submodule is also $\Z_2\,$-\,graded.
\newline\newline
If $\m_1,\m_2,\ldots,\m_n$ are \bk-supermodules then the direct
sum $\m_1\oplus\m_2\oplus\cdots\oplus\m_n$ is a \bk-supermodule
via the action
\[
    a(m_1,m_2,\cdots,m_n)=(am_1,am_2,\cdots
am_n)
\]
for all $a\in\bk,\, m_i\in\m_i,\, i=1,2,\cdots,n$.  A map
$\phi:\m_1\times\m_2\times\ldots\times\m_n \rightarrow\n$ where
\ns is a \bk-supermodule is said to have \emph{degree} $\alpha\
(\alpha\in\Z_2)$ if
\[
    \phi(m_1,\ldots,m_n)\in\n_{\alpha+|m_1|+\cdots+|m_n|}
\]
for all homogeneous $m_i\in\m_i,\,1\leq i\leq n$. Given
\bk-supermodules \ms and \n, we say that a map $\phi:\m\rightarrow
\n$ is \bk-linear if $\phi$ is additive and
\[
    \phi(ma)=\phi(m)a\textup{ for all
}a\in\bk,\,m\in\m.
\]
A \bk-linear map between \bk-supermodules \ms and \ns is called a
\emph{\bk-supermodule homomorphism from \ms to \n}.  The set of
all \bk-supermodule homomorphisms from \ms to \ns is denoted by
$\hm_{\bk}(\m,\n)$ and is a \bk-supermodule where the \bk-action
is given by $(af)(m):=a(f(m))$ and $(fa)(m):=f(am)$ for all
$a\in\bk,\,f\in\hm_{\bk}(\m,\n)$.  The $\Z_2$-grading is given by
\[
    \hm_{\bk}(\m,\n)=(\hm_{\bk}(\m,\n))_{\bz}\oplus(\hm_{\bk}(\m,\n))_{\bo}
\]
where $\hm_{\bk}(\m,\n)_{\alpha}=\{f\in\hm_{\bk}(\m,\n)\mid
\text{degree of }f\text{ is }\alpha\}$.  In case $\m=\n$ we write
$\ed_{\bk}\m$ instead of $\hm_{\bk}(\m,\m)$.  Given a
\bk-super\-module \p, we say that a map
$\phi:\m\times\n\rightarrow\p$ is \emph{\bk-bilinear} if
\begin{eqnarray}
    \phi(m+m',n+n')&=&\phi(m,n)+\phi(m,n')+\phi(m,n')+\phi(m',n');\\
    \phi(ma,n)&=&\phi(m,an);\text{ and}\\
    \phi(m,na)&=&\phi(m,n)a
\end{eqnarray}
for all $a\in \bk,\,m\in\m,\,n\in\n$.  The set of all \bk-bilinear
maps $\phi:\m\times\n\rightarrow\p$ forms a \bk-supermodule via
the action
\begin{eqnarray*}
    (a\phi)(m,n)&:=&a(\phi(m,n));\text{ and }\\
    (\phi a)(m,n)&:=&\phi(am,n)
\end{eqnarray*}
for all $a\in\bk,\,m\in\m,\,n\in\n$.  \newline\newline A
\bk-bilinear map $\phi:\m\times\m\rightarrow\n$ is
\emph{supersymmetric} if
\[
    \phi(m,n)=(-1)^{|m||n|}\phi(n,m)
\]
for all $m,n\in\m$.  Note that if $\phi$ is supersymmetric, then
$\phi|_{\m_{\bz}\times\m_{\bz}}$ is symmetric and
$\phi|_{\m_{\bo}\times\m_{\bo}}$ is skewsymmetric.  The
\emph{radical \calrs of} a supersymmetric $\phi$ is defined by
\[
    \calr=\{m\in\m \ve \phi(m,n)=0 \ \forall n\in\m\}.
\]
We say that $\phi$ is \emph{nondegenerate} if
$\calr=\{0\}$.\newline\newline A \emph{K-quadratic map} between
\bk-supermodules \ms and \ns is a \bk-supersymmetric, \bk-bilinear
map $q:\m\times\m\rightarrow\n$ of degree \bz.  Given a
\bk-supermodule \m, a \bk-quadratic map \qs from \ms to \bks is
called a \emph{K-quadratic form on \m}. Given two \bk-supermodules
$\m_1$ and $\m_2$ with \bk-quadratic forms $q_1$ and $q_2$
respectively, we define their \emph{orthogonal sum} $q=q_1\oplus
q_2$ to be the \bk-quadratic form on the \bk-supermodule
$\m=\m_1\oplus\m_2$ given by
\begin{eqnarray*}
    q(m_1\oplus m_2,n_1\oplus n_2)&=&(q_1\oplus q_2)(m_1\oplus m_2,n_1\oplus n_2)\\
    &=&q(m_1,n_1)+q(m_2,n_2)
\end{eqnarray*}
for $m_1,n_1\in\m_1,\,m_2,n_2\in\m_2$.
\begin{eg}\label{hyperbolic}\textup{Let $I$ be an arbitrary set.  We
define $H(I,K)$ to be the direct sum of free \bk-supermodules
$Kh_{\pm i}$ for $i\in I$ where $h_{\pm i}$ are of even degree.
In other words,
\[
    H(I)=H(I,K)=\oplus_{i\in I}(Kh_i\oplus Kh_{-i}).
\]
where the elements of the basis $\{h_{\pm i}\}_{i\in I}$ are all
of even degree. Then we have that $H(I)=H(I)_{\bz}\oplus
H(I)_{\bo}$ is a \bk-supermodule where
\[
    H(I)_{\alpha}=\oplus_{i\in I}(K_{\alpha}h_i\oplus K_{\alpha}h_{-
i})
\]
for $\alpha\in\Z_2$.  We call $H(I,K)$ the \emph{hyperbolic
\bk-superspace}.  In the case when $|I|=1$, we call $H(I,K)$ the
\emph{hyperbolic K-superplane}.  We define the \bk-quadratic form
$q_{I}$ on the hyperbolic \bk-superspace $H(I)$ by setting
\[
    q_{I}(h_{\sigma i},h_{-\mu j})=\delta_{\sigma,\mu}\delta_{i,j}
\]
for all $i,j\in I,\,\sigma,\mu\in\{\pm\}$ and extending $\q_{I}$
bilinearly over \bk.  Note that $q_{I}$ is the orthogonal sum of
all $q_{\{i\}},\ i\in I$. }\end{eg} A \emph{superalgebra \bas over
\bk} is a \bk-supermodule with a \bk-bilinear map $\cdot:A\times
A\rightarrow A$ of degree \bz.  A \emph{subalgebra B of the
superalgebra A} is a \bk-submodule $B$ of $A$ such that $B$ is
closed under $\cdot|_{B\times B}$.  A subalgebra $B$ of $A$ is an
\emph{ideal} of $A$ if $A\cdot B\subset B$.  A unital,
associative, supercommutative superalgebra $A$ over \bks is called
a \emph{superextension of K}.\newline\newline
\begin{eg}\label{Grassmann}\textup{
Let $\g$ be the associative \bk-algebra generated by
$\{\xi_i\}_{i\in\N}$ subject to the relation
$\xi_i\xi_j+\xi_j\xi_i=0\text{ for all } i,j\in\N$. Since
$\half\in\bk$, we have in particular that $\xi_i^2=0$ for all
$i\in\N$. The algebra $\g$ is called the \emph{exterior algebra
over \bks on a countable number of generators} $\xi_i,\,i\in\N$,
or simply the \emph{Grassmann algebra}.  It is easy to check that
the Grassmann algebra is a superextension of \bk. }\end{eg} Let
\ms be a \bk-supermodule.  Then $\ed_K\m$ is a
\bk-superalgebra via the usual composition of maps.\\
\\
Given two \bk-superalgebras $A$ and $B$, a \bk-algebra
homomorphism $\phi:A\rightarrow B$ of degree \bz\ is called a
\emph{\bk-superalgebra homomorphism}.  If, in addition, $\phi$ is
bijective, we say that $\phi$ is an \emph{isomorphism} and we
write $A\overset{\phi}{\cong}B$, or simply $A\cong B$.
\begin{Defn}\label{module}\textup{
Let \bas be a \bk-superalgebra in a variety $\mathfrak{V}$ (later
$\mathfrak{V}$ will be the variety of Jordan or Lie
superalgebras).  Let \cus be a \bk-supermodule equipped with a
pair of \bk-bilinear mappings $(a,u)\mapsto au,\, (a,u)\mapsto
ua,\,a\in\ba,\,u\in\cu$, of $\ba\times\cu$ into \cus of degree 0.
Then $X=\ba\oplus\cu$ is a \bk-supermodule on which we define a
multiplication by
\begin{eqnarray}
    (a+u)(b+v)=ab+av+ub
\end{eqnarray}
for all $a,b\in\ba,\,u,v\in\cu$.  Since this product is
\bk-bilinear, $X$ is a superalgebra over \bk, called the
\emph{split null extension} of \cas determined by the bilinear
mappings of \bas and \cus (see \cite[Chap.\,II,\ Sect.,\,5]{J} for
the classical case).  If $X$ is a superalgebra in the variety
$\mathfrak{V}$ then we say that $\cu$ is a
\emph{$\mathfrak{V}$-module} for \ba.  In this case, for each
$\alpha\in\Z_2$ let $(\der_K(\ba,\cu))_{\alpha}$ be the space of
all homogeneous \bk-module homomorphisms
$d\in\hm_{\bk}(\ba,\cu)_{\alpha}$ satisfying for all homogeneous
$x\in\ba$ and all $y\in \ba$
\[
    d(xy)=d(x)y+(-1)^{\alpha|x|}xd(y).
\]
We define
\[
    \der_\bk(\ba,\cu)=(\der_{\bk}(\ba,\cu))_{\bz}\oplus(\der_{\bk}(\ba,\cu))_{\bo}.
\]
Then it is easy to see that $\der_{\bk}(\ba,\cu)$ is a submodule
of $\ed_{\bk}\m$ and hence a \bk-supermodule.  The elements of
$\der_{\bk}(\ba,\cu)$ are called the \emph{\bk-derivations (from
\bas to \cu)}.  If $\ba=\cu$, we simply write $\der_{\bk}\ba$.
}\end{Defn} A \bk-supermodule $L$ is called a \emph{Lie
superalgebra} if it is equipped with a \bk-bilinear bracket
multiplication $[\cdot,\cdot] :L\times L\rightarrow L$ such that:
\begin{description}
\item[(SL1)]  For all $x,y\in L$,
\[
    [x,y]=-(-1)^{|x||y|}[y,x];
\]
\item[(SL2)]  For all $x,y,z\in L$,
\[
    (-1)^{|x||z|}[[x,y],z]+(-1)^{|y||x|}[[y,z],x]+(-
1)^{|z||y|}[[z,x],y]=0.
\]
\end{description}
The property (SL2) is called the \emph{Jacobi identity} and is
equivalent to
\begin{eqnarray}
    \label{Jacobi}[[x,y],z]=[x,[y,z]]-(-1)^{|x||y|}[y,[x,z]]
\end{eqnarray}
for all $x,y,z\in L$.\newline\newline For example, given an
associative \bk-superalgebra $A$, we can define a Lie superalgebra
structure on $A$, denoted by $A^{(-)}$, by defining the bracket
operation via $[x,y]=xy-(-1)^{|x||y|}yx$ for all $x,y\in A$.  In
particular, one can easily check that $\der_KA$ is a subalgebra of
$(\ed_K A)^{(-)}$.\newline\newline If $L$ is a Lie superalgebra,
we define the \emph{centre of L} by
\[
    Z(L)=\{x\in L\mid [x,y]=0\text{ for all }y\in L\}.
\]
We have that $Z(L)$ is an ideal of $L$.  A Lie superalgebra $L$ is
called a \emph{central extension} of the Lie superalgebra
$\tilde{L}$ if there exists a central ideal $C$ of $L$ (i.e.,
$C\subset Z(L)$) such that $L/C\cong\tilde{L}$.\newline\newline
For any element $x$ of $L$, we can define the map
$\ad(x):L\rightarrow L$ by $\ad(x)(y):=[x,y]$.  The
\bk-supermodule $\ad(L)$ forms an ideal of $\der_K L$ called the
\emph{inner derivation algebra of L} and denoted by $\ider L$. For
example, any finite-dimensional Lie superalgebra $L$ over a field
\lk, whose Killing form is nondegenerate, has no outer
derivations, i.e.\! $\der_{\lk}L=\ider L$ (see
\cite[Prop.\,2.3.4]{Kac2}).\newline\newline We recall the
following definition (see \cite{N1}, for an equivalent definition
see \cite{LN}).  A subset $R$ of a real vector space $X$ with a
scalar product $(\cdot,\cdot)$ is called a \emph{root system (in
X)} if $R$ has the following three properties:
\begin{description}
\item[(RS1)]  $R$ generates $X$ as a vector space and $0\notin R$,
\item[(RS2)]  for each $\alpha\in R$ we have $s_{\alpha}(R)=R$
where $s_{\alpha}(x)=x-2\frac{(x,\alpha)}{(\alpha,\alpha)}\alpha$
for $x\in X$, \item[(RS3)]  for all $\alpha,\beta\in R$, we have
$\langle\alpha,\beta\rangle\in\Z$, where for $x\in X,\ \alpha\in
R$,
\[
    \langle x,\alpha\rangle:=2\frac{(x,\alpha)}{(\alpha,\alpha)}.
\]
\end{description}
A root system $R$ is \emph{irreducible} if $R\neq 0$ and if $R$ is
not an orthogonal sum of two non-empty root systems.  Every root
system is an orthogonal sum of irreducible root systems which are
uniquely determined and are called the \emph{irreducible
components of R}.  Note that in the finite-dimensional case, these
definitions conform with the traditional definition of (finite)
root systems.  An irreducible root system is isomorphic to a
finite root system or to one of the infinite analogues of the root
systems of type $A-D$ and $BC$ (\cite{LN}).\newline\newline Let
$R$ be a root system.  We say that a Lie superalgebra over \bks is
\emph{$R$-graded} (\cite{GN}) if there exist supermodules
$L^{\alpha},\ \alpha\in R\cup\{0\}$ of $L$, such that
\begin{itemize}
    \item[(SG1)]  $L=\oplus_{\alpha\in R\cup\{0\}}L^{\alpha}$;
    \item[(SG2)]  for all $\alpha,\beta\in R\cup\{0\}$,
\[
   [L^{\alpha},L^{\beta}]\subset\left\{\begin{array}{ll}L^{\alpha+\beta}&\textup{if\ }
\alpha+\beta\in R\cup\{0\} \\
                \{0\}&\textup{if\ }\alpha+\beta
\notin R\cup\{0\}
            \end{array}\right.;
\]
\item[(SG3)]  as a Lie superalgebra, $L$ is generated by
$\cup_{\alpha\in R}L^{\alpha}$; \item[(SG4)]  for every $\alpha\in
R$ there exists $0\neq X_{\alpha}\in (L^{\alpha})_{\bz}$ such that
$H_{\alpha}:=[X_{-\alpha},X_{\alpha}]$ operates on $L^{\beta}\
(\beta\in R\cup\{0\})$ by
\[
    [H_{\alpha},z_{\beta}]=\langle\beta,\alpha\rangle z_{\beta}\quad
(z_{\beta}\in L^{\beta})
\]
where
\[
    \langle\beta,\alpha\rangle=2\frac{(\beta,\alpha)}{(\alpha,\alpha)}.
\]
\end{itemize}
This is a generalization of $ R $-graded Lie algebras (see
\cite{BM} and \cite{BeZ}) in the following sense:  Every
$R$-graded Lie algebra is an $R$-graded Lie superalgebra (by
setting $L_{\bz}= L$ and $L_{\bo}={0}$).  However, if $L$ is an $
R $-graded Lie superalgebra, $L_{\bz}$ need not be an $ R $-graded
Lie algebra (see Example \ref{countereg}).
\newline\newline
A \emph{Jordan superalgebra} (\cite{Kac3}) $J=J_{\bz}\oplus
J_{\bo}$ over \bks is a \bk-superalgebra such that for all
$a,b,c,d\in J$,
\begin{description}
\item[(JSA1)]  $ab=(-1)^{|a||b|}ba$ \item[(JSA2)]
$((ab)c)d+(-1)^{|d||c|}a((bd)c)+(- 1)^{|b||a|+|d||c|}b((ad)c)= $
\item[\qquad\qquad\qquad ]  $(ab)(cd)+(-1)^{|c||b|}(ac)(bd)+(-
1)^{|d|(|b|+|c|)}(ad)(bc)$ \item[(JSA3)]  $(a^2c)a=a^2(ca)$.
\end{description}
\begin{eg}\label{JSA of q.f.}\textup{Let \ms be a \bk-supermodule
with a \bk-quadratic form \q.  Then the algebra $J=\bk\oplus\m$
with multiplication $\circ$ given by
\[
    (a\oplus m)\circ(b\oplus n)=(ab+q(m,n))\oplus(an+mb)
\]
forms a Jordan superalgebra, called the Jordan superalgebra
associated with the quadratic form \q. }\end{eg} One can check
that given a Jordan algebra $J$ over \bk, a \bk-module \ms is a
Jordan superalgebra module (or a \emph{Jordan supermodule})
according to Definition \ref{module} if \ms is a \bk-supermodule
satisfying:
\begin{description}
\item[(JSAM1)]  $am=(-1)^{|a||m|}ma$ \item[(JSAM2)]
$0=(-1)^{|b||c|}[R_{ab},R_c]+(-1)^{|a||b|}[R_{ca},R_b]+(-
1)^{|a||c|}[R_{bc},R_a]$ \item[(JSAM3)]
$R_cR_bR_a+(-1)^{|a||b|+|a||c|+|b||c|}R_aR_bR_c+(-
1)^{|b||c|}R_{(ac)b}=$ \item[\ ]
\qquad\qquad\qquad$(-1)^{(|a|+|b|)|c|}R_{ab}R_c+(-
1)^{|a||b|}R_{ac}R_b+(-1)^{|b||c|}R_{cb}R_a$
\end{description}
for all $a,b,c\in J,\, m\in\m$ where $R_d$ denotes the action on
the right of $J$ on $X$, i.e., $R_dx=xd$.


\section{ORTHOSYMPLECTIC LIE SUPER\-ALGEBRAS}
We will denote by \bas a superextension of \bks and by \ms an
\ba-supermodule with an \ba-quadratic form \qs.  \newline\newline
Set
\begin{align*}
    \osp(q)=\{x\in\ed_{\ba}\m\mid &q(x(m),n)+(-1)^{|m||n|}q(x(n),m)=0
\text{ for all }\\
    &m,n\in\m\}.
\end{align*}
We define, for homogeneous $m,n,p\in\m$,
\begin{eqnarray}\label{defn of E}
    \e_{m,n}(p)&=&mq(n,p)-(-1)^{|n||p|}q(m,p)n.
\end{eqnarray}
and we extend $\e:\m\times\m\rightarrow\ed_{\ba}\m$ linearly.  Let
\[
    \eosp(\q)=\text{span}_{\Z}\{\e_{m,n}\mid m,n\in\m\}.
\]
We call \osp(\q) the \emph{orthosymplectic Lie superalgebra (of
q)} and \eosp(\q) the \emph{elementary orthosymplectic Lie
superalgebra (of q)}. \newline\newline The techniques developed in
\cite{N2} to describe root-graded Lie algebras can be adapted to
the supercase, see \cite[Section 5.4]{N2} and \cite{GN}.  In
particular, let $I$ be an arbitrary set and $q_I$ the quadratic
form on $H(I;K)$ defined as in Example \ref{hyperbolic}.
\begin{description}
\item[($B_I$)]  Assume that $\frac{1}{6}\in\bk$ and $|I|\geq 3$.
A Lie superalgebra $L$ over \bks is $B_I$-graded if and only if
there exists a superextension \bas of \bk, an \ba-supermodule \ms
with an \ba-quadratic form \qs and an element $m\in\m$ satisfying
$q(m,m)=1$ such that $L$ is a central extension of the Lie
superalgebra $\eosp(q_I\oplus q)$. \item[($D_I$)]  Assume that
$\frac{1}{6}\in\bk$ and $|I|\geq 4$.  A Lie superalgebra $L$ over
\bks is $D_I$-graded if and only if there exists a superextension
\bas of \bks such that $L$ is a central extension of the Lie
superalgebra $\eosp(q_I)$.
\end{description}
\begin{eg}\label{countereg}\textup{Let \bks be a commutative and unital
ring containing $\frac{1}{6}$.  We consider \bks as a superring
with $\bk_{\bo}=(0)$. Let $I$ and $J$ be arbitrary index sets
($|I|\geq 3$ and $J\neq\emptyset$) and let $\m=\oplus_{j\in J}\bk
m_j$ be a free \bk-supermodule with a basis consisting of odd
elements and with a nonzero \bk-quadratic form $q_{\m}$.  Consider
$L=\eosp(q_I\oplus q_{\m}\oplus q_0)$ where $q_0$ is the quadratic
form on the free module $\bk x_0$ given by $q_0(x_0,x_0)=1$. Then
by \emph{($B_I$)}, $L$ is a $B_I$-graded Lie superalgebra over
\bk. It can be shown in particular that
$\eosp(q_{\m})(\hookrightarrow\eosp(q_I\oplus q_{\m}\oplus q_0)$)
is contained in $(L^0)_{\bz}$ but that it is not generated by
$(L^{\alpha})_{\bz},\ \alpha\in B_I$.  Hence $L_{\bz}$ is not a
$B_I$-graded Lie algebra. }\end{eg} If $\und{x}=(x_i)_{i\in I}$ is
a family in $\eosp(q)$ we say that $\und{x}$ is \emph{summable} if
for every $m\in\m,\ (x_i(m))_{i\in I}$ has only finitely many
nonzero terms, and in this case we define $x=\sum_{i\in
I}x_i\in\ed_{\ba}\m$ by setting
\[
    \left(\sum_{i\in I}x_i\right)(m)=\sum_{i\in I_m}x_i(m)
\]
for $I_m=\{i\in I\mid x_i(m)\neq 0\}$ ($I_m$ is finite).  We will
say that $\sum_{i\in I}x_i$ is summable to indicate that
$(x_i)_{i\in I}$ is summable.  Also, given a summable
$x=\sum_{i\in I}x_i$ we will throughout denote by $I_m$ the finite
subset of $I$ such that for $i\in I,\ x_i(m)\neq 0$.  The set of
all summable $\sum_{i\in I}x_i,\ x_i\in\eosp(q)$ is denoted by
$\ovl{\eosp(q)}$ and is an \ba-submodule of $\ed_{\ba}\m$ with the
natural $\Z_2$-grading. Summable homomorphisms are discussed in a
more general setting in \cite{D}.  In this section we will
investigate some of the properties and the structures of
\eosp(\q), $\ovl{\eosp(q)}$ and \osp(\q).
\begin{lemma}  \label{properties of E}Let \ms be an \ba-supermodule with
an \ba-quadratic form \q. Then
\begin{itemize}
\item[1)]  for all $m,n\in\m$ and for all $a\in\ba$,
$\e_{m,n}\in\osp(q)$.  Moreover, the map
$\e:\m\times\m\rightarrow\ed_A\m:(m,n)\mapsto\e_{m,n}$ is
\ba-bilinear,
\[
    \e_{m,n}=-(-1)^{|m||n|}\e_{n,m}; \text{\ and}
\]
\[
    q(\e_{m,n}(p),r)=(-1)^{(|m|+|n|)(|p|+|r|)}q(\e_{p,r}(m),n)
\]
for all $m,n,p,r\in\m$; \item[2)]  \osp(\q) is a subalgebra of the
Lie superalgebra $(\ed_{\ba}\m)^{(-)}$; \item[3)]  for all
$m,n\in\m$ and $x\in\osp(q)$,
\begin{eqnarray}\label{[osp,E]}
    [x,\e_{m,n}]&=&\e_{x(m),n}+(-1)^{|x||m|}\e_{m,x(n)}; \ and\
hence
\end{eqnarray}
\item[4)]  $\eosp(\q)$ and $\ovl{\eosp(q)}$ are ideals of
$\osp(\q)$. \item[5)]  If \bas is a field of characteristic 0, \qs
is nondegenerate and \ms is a free \ba-module of finite rank, then
$\osp(q)=\ovl{\eosp(q)}=\eosp(q)$.
\end{itemize}
\end{lemma}
The easy proof of this lemma is left to the reader.  However, we
note that to prove that $\ovl{\eosp(q)}$ is an ideal of \osp(\q)
one uses the formula
\begin{eqnarray}\label{ovl ideal}
    \left[\sum_{i\in I}x_i,y\right]=\sum_{i\in I}[x_i,y]\in\ovl{\eosp(q)},\qquad
\sum_{i\in I}x_i\in\ovl{\eosp(q)},\,y\in X
\end{eqnarray}
for any subalgebra $X$ of $\ed_{\bk}\m$ of which $\eosp(q)$ is an ideal.\\
\\
In what follows we will denote by \ces a subalgebra of \osp(\q)
containing \eosp(\q).  We set
\[
    \ce_{\infty}=\{(a,m,n,x)\mid a\in\ba,\ m,n\in\m,\ x\in\ce\}
\]
to be the Lie superalgebra with multiplication
\begin{align}
\nonumber   [(&a,m,n,x),(a',m',n',x')]=(q(m,n')-q(n,m'),\\
\nonumber   &am'+xm'-ma'-(-1)^{|x'||m|}x'm,-an'+xn'+na'-(-1)^{|x'||n|}x'n,\\
\label{mult}    &\e_{m,n'}+\e_{n,m'}+[x,x'])
\end{align}
for all $a,a'\in \ba,\ m,m',n,n'\in\m,\ x,x'\in\ce$. Then
$L=\ce_{\infty}=L_+\oplus L_0\oplus L_-$ is a 3-graded Lie
superalgebra (i.e., $L_{\sigma}L_{\mu}\subset L_{\sigma+\mu},\
\sigma,\mu\in\{\pm,0\}$) where
\[
\begin{array}{rcl}
L_0&=&\{(a,0,0,x)\mid a\in \ba,\ x\in\eosp(q)\}\\
L_+&=&\{(0,m,0,0)\mid m\in\m\}\\
L_-&=&\{(0,0,m,0)\mid m\in\m\}.
\end{array}
\]
If we set $\q_{\infty}:=\q_{\{\infty\}}\oplus\q$ where
$q_{\{\infty\}}$ is the \ba-quadratic form on the hyperbolic
superplane $H(\{\infty\},\ba)$ then
$\eosp(\q)_{\infty}=\eosp(\q_{\infty})$,
$\osp(q)_{\infty}=\osp(q_{\infty})$ and so
$\eosp(q_{\infty})\subset\ce_{\infty}\subset\osp(q_{\infty})$. In
particular, we see from ($B_I$) and ($D_I$) that $B_I$-graded and
$D_I$-graded Lie superalgebras ($|I|\geq 3$ and $|I|\geq 4$
respectively) over \bks containing $\frac{1}{6}$ are central
extensions of $\eosp(q)_{\infty}$ for a suitable superextension
\bas of \bk, an \ba-supermodule \ms and an \ba-quadratic form \q.
\newline\newline Let \ms be a \bk-supermodule with a \bk-quadratic
form \q.  Let $J=\bk\oplus\m$ be the Jordan superalgebra
associated with the quadratic form \qs (see Example \ref{JSA of
q.f.}). Consider the action of $J$ on the \bk-supermodule
$X=\ce\oplus \m$ given by
\begin{eqnarray}
\label{jsm} (x\oplus p)(a\oplus m)&=&(xa+\e_{p,m})\oplus(x(m)+pa)\\
\nonumber   &=&(-1)^{|a\oplus m||x\oplus p|}(a\oplus m)(x\oplus p)
\end{eqnarray}
for $a\in\bk,\,m,p\in\m$ and $x\in\ce$.  Then one can easily check
that (JSAM1)-(JSAM3) hold and so $X$ is a Jordan superalgebra
module for $J$.  \newline\newline The following can be easily
shown and will be used in the next section:
\begin{prop}
\label{centre} Let \ms be an \ba-supermodule with an \ba-quadratic
form \q.  Let $L=\ce_{\infty}$.  Then $Z(L)=(\ann_A\m,0,0,0)$
where $Ann_A\m=\{a\in A\mid am=0 \textup{\ for all\ }m\in\m\}$.
\end{prop}


\section{DERIVATIONS}
In this section we will study the derivations of $\ce_{\infty}$
where \ces is a subalgebra of \osp(\q) containing \eosp(\q) for an
\ba-supermodule \ms with an \ba-quadratic form \q.  We will denote
$\ce_{\infty}$ by $L$ throughout.  Recall that $L=\{(a,m,n,x)\mid
a\in A,\,m,n\in\m,\,x\in\ce\}$ with multiplication given by
(\ref{mult}).
\begin{prop}
\label{der for toral L0=A+[L+,L-]} We have that
\[
    \der_{\bk} L=\ider L+\{d\in\der_{\bk} L\mid d(1,0,0,0)\in Z(L)\}.
\]
Moreover, \[\begin{array}{lrl}
    (\der_{\bk} L)_0&:=&\{d\in\der_KL\mid d(L_{\sigma})\in
L_{\sigma},\sigma=0,\pm\}\\
    &=&\{d\in\der_{\bk} L\mid d(1,0,0,0)\in Z(L)\}.
    \end{array}\]
    In particular, $\{d\in\der_{\bk} L\mid d(1,0,0,0)\in Z(L)\}$ is a subalgebra of
$\der_KL$ and
\[
    \ider L\cap\{d\in\der_KL\mid d(1,0,0,0)=0\}= \ad L_0.
\]
\end{prop}
\textsl{Proof:} It is clear that $\ider L+(\der_{\bk}
L)_0\subset\der_{\bk} L$.  For $x\in L$ we will write
$x=x_-+x_0+x_+$ where $x_{\mu}\in L_{\mu},\ \mu=0,\pm$ and we will
denote $(1,0,0,0)$ simply by 1.  Let $d\in\der_{\bk}L$ and set
$d'=d+\ad(d(1)_+-d(1)_-)$.  Then for all $x\in L$
\begin{eqnarray*}
    d'(x_+)&=&d(x_+)+[d(1)_+-d(1)_-,x_+]=d[1,x_+]-
[d(1)_-,x_+]\\
    &=&[d(1),x_+]+[1,d(x_+)]-[d(1)_-,x_+]\\
    &=&\underbrace{[(d(1))_0,x_+]}_{\in
L_+}+\underbrace{[1,d(x_+)_+]}_{\in L_+}\in L_+
\end{eqnarray*}
since
\[
    d(x_+)_-=(d[1,x_+])_-=\underbrace{[d(1),x_+]_-
}_{=0}+[1,d(x_+)]_-=-d(x_+)_-
\]
and so $d(x_+)_-=0$.  Similarly, $d'(x_-)\in L_-$. In addition,
$[1, x_0]=0$ implies $0=[d(1),x_0]+[1,d(x_0)]$ which in turn
implies $d(x_0)_{\pm}=\mp[d(1)_{\pm},x_0]$ and hence
$d'(x_0)=d(x_0)_0\in L_0$.  Therefore $\der_{\bk}L=\ider
L+(\der_{\bk} L)_0$.  Now, let $d\in(\der_{\bk}L)_0$ and $x\in L$.
Then
\[
    d(x_+)=d[1,x_+]=[d(1),x_+]+[1,\underbrace{d(x_+)}_{\in
L_+}]=[d(1),x_+]+d(x_+)
\]
which implies that $[d(1),x_+]=0$.  Similarly, $[d(1),x_-]=0$.
Moreover, since $[1, x_0]=0$, we have $0=[d(1),x_0]+[1,{d(x_0)}]$
and so since ${d(x_0)}{\in L_0},\ [d(1),x_0]=0$.  Hence
$(\der_{\bk}L)_0\subset\{d\in\der_KL\mid d(1)\in Z(L)\}$.
Conversely, if $d\in\der_{\bk}L$ is such that $d(1)\in Z(L)$ then
for all $x\in L$, $d(x_{\pm})=\pm d[1,x_{\pm}]$ implies
$d(L_{\pm})\subset L_{\pm}$ and $[1,x_0]=0$ implies $[1,d(x_0)]=-
[d(1),x_0]=0$ which in turn implies $d(L_0)\subset L_0$. Therefore
$(\der_{\bk}L)_0=\{d\in\der_KL\mid d(1)\in Z(L)\}$. It is easily
verified that $(\der_{\bk} L)_0$ is a subalgebra of $\der_KL$ and
so all of the statements of the proposition hold.
\qed\newline\newline We remark that the proposition above holds
for so called \emph{toral 3-graded Lie superalgebras}, i.e., for
the Lie superalgebras of the form $L=L_{-1}\oplus L_0\oplus L_1$
with $L_{\sigma}L_{\mu}\subset L_{\sigma+\mu}$ for
$\sigma,\mu\in\{0,\pm1\}$ and such that there exists $\zeta\in
L_0$ satisfying $[\zeta,x_{\sigma}]=\sigma x_{\sigma}$ for all
$x_{\sigma}\in L_{\sigma},\ \sigma\in\{0,\pm\}$.\newline\newline
In what follows, we will assume that $\ann_A\m=0$ and so by
Proposition \ref{centre}, $Z(L)=0$.  In addition we have the
superalgebra monomorphism $A\rightarrow \ed_A{\m}:a\mapsto
a\id|_{\m}$.\\
\\
We make the following definitions which will be instrumental in
describing the derivations of $\ce_{\infty}$.
\begin{Defn}\textup{
For each $\alpha\in\Z_2$, let $(\es_{\ce})_{\alpha}$ be the set of
all \bk-linear endomorphisms $S$ of $\m$ of degree $\alpha$
satisfying
\begin{itemize}
\item[\textbf{(S1)}]  $[S,\ba\id ]\subset\ba \id $;
\item[\textbf{(S2)}]  $[S,\q(m,n)\id ]=(\q(S(m),n)+(-
1)^{|m||n|}q(S(n),m))\id $ for all homogeneous $m,n\in \m$;
\item[\textbf{(S3)}]  $[S,\ce]\subset\ce$.
\end{itemize}
}\end{Defn}
\begin{Defn}\textup{
For each $\alpha\in\Z_2$, let $(\te_{\ce})_{\alpha}$ be the set of
all \bk-linear endomorphisms $T$ of $\m$ of degree $\alpha$
satisfying
\begin{itemize}
\item[\textbf{(T1)}]  $[T,\ba\id ]\subset \ce$;
\item[\textbf{(T2)}]  $[T,\q(m,n)\id ]=\e_{T(m),n}+(-
1)^{|m||n|}\e_{T(n),m}$ for all homogeneous $m,n\in \m$;
\item[\textbf{(T3)}]  $[T,\ce]\subset\ba\id$.
\end{itemize}
}\end{Defn} Set
$\es_{\ce}=(\es_{\ce})_{\bz}\oplus(\es_{\ce})_{\bo}\textup{ and }
\te_{\ce}=(\te_{\ce})_{\bz}\oplus(\te_{\ce})_{\bo}.$ We see
immediately that $\es_{\ce}$ and $\te_{\ce}$ are submodules of the
\bk-supermodule $\ed_K(\m)$.  Moreover, $\ce\subset \es_{\ce}$ and
$\ba\id\subset\te_{\ce}$. We will write \ess and \tes for
$\es_{\ce}$ and $\te_{\ce}$ respectively when the context is
clear.  In fact, if we define, for $N\in\{1,2,3\}$ and
$\alpha\in\Z_2$,
\[
    \es^{(N)}_{\alpha}=\{S\in(\ed_{\bk}\m)_{\alpha}\mid S\text{
satisfies (SN)}\}
\]
and
\[
    \te^{(N)}_{\alpha}=\{T\in(\ed_{\bk}\m)_{\alpha}\mid T\text{
satisfies (TN)}\}
\]
and set $\es^{(N)}=\es^{(N)}_{\bz}\oplus\es^{(N)}_{\bo},\
\te^{(N)}=\te^{(N)}_{\bz}\oplus\te^{(N)}_{\bo}$ then
\[
    \es_{\ce}=\bigcap_{N=1}^3\es^{(N)}\textup{\quad
and\quad}\te_{\ce}=\bigcap_{N=1}^3\te^{(N)}.
\]
\begin{rmk}\textup{
Note that the $(\ )^{(N)}$ notation is traditionally used to
denote the derived algebras of a Lie algebra.  However in what
follows we will not be using derived algebras at all so there
should be no confusion resulting from the above definition.
}\end{rmk} We have the following properties of \ess and \te:
\begin{prop}  \label{prop of S and T}\
\begin{itemize}
\item[(1)]  Properties of $\es_{\ce}$:
\begin{itemize}
\item[(i)]  For all $S\in\esone$, the map $\Delta:A\rightarrow A$
given by $\Delta(a)\id=[S,a\id]$ is a \bk-derivation of \ba.
Moreover, $[\ed_A\m,S]\subset\ed_A\m$. \item[(ii)]  For all
$S\in\estwo$ and for all $m,n\in\m$,
\begin{eqnarray}
\label{[S,E]}   [S,\e_{m,n}]=\e_{S(m),n}-(-1)^{|m||n|}\e_{S(n),m}.
\end{eqnarray}
\item[(iii)]   Let $X=\eosp(\q)$, $\ovl{\eosp(q)}$ or \osp(\q).
Then $[\estwo,X]\subset X$.  In particular, if $\ce=\eosp(\q)$,
$\ovl{\eosp(q)}$ or \osp(\q) then $\estwo\subset\esthree$.
\item[(iv)] $\osp(q)$ is an ideal of $\esone\cap\estwo$ and so
\begin{eqnarray}
\label{osp in S}    \osp(q)\cap\esthree\lhd\es_{\ce}.
\end{eqnarray}
In particular, $\ce\lhd\es_{\ce}$.  Moreover,
$\osp(q)\lhd\es_{\ce}$ for any ideal \ces of \osp(\q) which is
true in particular for $\ce=\eosp(\q),\ \ovl{\eosp(q)}$ and
\osp(\q).
\end{itemize}
\item[(2)] Properties of $\te_{\ce}$:
\begin{itemize}
\item[(i)]  For all $T\in\teone$, the map
$\varphi:A\rightarrow\ce$ given by $\varphi(a)=[T,a\id]$ is a
\bk-derivation from \bas to \ce. \item[(ii)]  For all $T\in\tetwo$
and for all $m,n\in\m$,
\begin{eqnarray}
\label{[T,E]}
\qquad[T,\e_{m,n}]=(q(T(m),n)-(-1)^{|m||n|}q(T(n),m))\id.
\end{eqnarray}
\item[(iii)]  Suppose that $\ce\subset\ovl{\eosp(q)}$ and
$\ann_{\ba}(m)= 0$ for some $m\in\m$ if $\ce\neq\eosp(q)$.  Then
$\tetwo\subset\tethree$. \item[(iv)]
$[\ce,[\tethree,A\cdot\idd]]=(0)$ and $[\tethree,[\ce,\ce]]=(0)$.
\end{itemize}
\end{itemize}
\end{prop}
\begin{rmk}\textup{
We should note that the assertion (1iii) in the above proposition
does not hold for general $\eosp(q)\subset\ce\subset\osp(q)$, for
an example see Corollary \ref{der of Kac}.  Also, a slightly more
general version of (2iii) is proved in \cite{D}. }\end{rmk}
\textsl{Proof:}  We leave it to the reader to check the above
statements.  They follow from the properties \sone-\sthrees and
\tone-\tthree, except for the following two cases:  in (1iii) one
needs to verify that $[\estwo, \eosp(q)]\subset\eosp(q)$ implies
$[\estwo, \ovl{\eosp(q)}]\subset\ovl{\eosp(q)}$ for which one uses
(\ref{ovl ideal}), and in (2iii) one needs to check that for $T\in
\tetwo$ and $x\in\ovl{\eosp(q)}$, $[T,x]=\sum_i a_i\id$ for some
index set $I$ where for each $n\in\m$ there are only finitely many
$i\in I$ such that $a_in\neq 0$; this, together with the fact that
$\ann_{\ba}(m)=0$ for some $m\in\m$, implies that in fact only
finitely many $a_i,\,i\in I$ are nonzero and so
$[\tetwo,\ovl{\eosp(q)}]\subset \ba\id$.\qed\newline
\newline
In fact, \ess and \tes give rise to a Lie superalgebra structure
in the following way:
\begin{prop}\label{S+T is LSA}  \
The \bk-supermodule $\es\oplus\te\subset \ed_K\m\oplus\ed_K\m$ is
a Lie superalgebra with multiplication
\begin{eqnarray}
\label{mult in S+T} [S\oplus T,S'\oplus
T']\til=([S,S']+[T,T']])\oplus([S,T']+[T,S'])
\end{eqnarray}
where $[\ ,\ ]$ is the usual bracket multiplication in
$(\ed_K\m)^{(- )}$.
\end{prop}
\textsl{Proof:}  By Proposition \ref{prop of S and T}, the bracket
$[\ ,\ ]\til$ is well-defined.  It is easy to check that
$\es\oplus\te=(\es_{\bz}\oplus\te_{\bz})\oplus(\es_{\bo}\oplus\te_{\bo})$
defines a $\Z_2$-grading which gives $\es\oplus\te$ a
\bk-superalgebra structure under the multiplication given by
(\ref{mult in S+T}).  The fact that $(\ed_K\m)^{(-)}$ is a Lie
superalgebra now implies the assertion.\qed
\begin{rmk}\textup{It is important to note that $\es\oplus\te$ does not
in general naturally imbed in $\ed_K\m$.  For example, if
$\q\equiv 0$ then $\ed_A\m\subset\es\cap\te$ and hence
$\es\cap\te\neq0$.}
\end{rmk}
\begin{Defn}\textup{
Let $\de_{\m}$ be the set of all \bk-endomorphisms $d_{S,T}$ of
$L$ for which there exist $S\in\es,\ T\in\te$ satisfying
\begin{eqnarray}
    \label{dformula}d_{S,T}(a,m,n,x)&=&([S,a\id]+[T,x],(S+T)(m),\\
    &&(S-T)(n),[T,a\id]+[S,x]).\nonumber
\end{eqnarray}
}\end{Defn}
\begin{rmk}\label{uniqueness in D}\textup{
Note that for every $d\in\de_{\m}$, $d=d_{S,T}$ for unique
$S\in\es,\ T\in\te$. }\end{rmk} We will often use the following
notation:
\begin{eqnarray}\label{abbrev}\begin{array}{lcll}
    a&:=&(a,0,0,0)&\text{for all }a\in\ba;\\
    m^+&:=&(0,m,0,0)&\text{for all }m\in\m;\\
    m^-&:=&(0,0,m,0)&\text{for all }m\in\m;\\
    x&:=&(0,0,0,x)&\text{for all }x\in\eosp(q);
\end{array}\end{eqnarray}
the meaning will be clear from the context.
\begin{thm}\label{main}\ \newline
1)  We have $\de_{\m}=(\der_KL)_0$.  In particular, $\de_{\m}$ is
a subalgebra of $\der_K L$ and the map
$\phi:\es\oplus\te\rightarrow\de_{\m}:S\oplus T\mapsto d_{S,T}$ is
a Lie superalgebra isomorphism.\newline 2)  Moreover,
$\der_KL=\ider L+\de_{\m}$ and\newline 3)  we have that $\ider
L\cap\de_{\m}=\{d_{S,T}\in\de_{\m}\mid S\in\ce,\ T\in A\id\}=\ad
L_0$.\newline 4)  If $\es=\ce\oplus\es_0$ and
$\te=A\id\oplus\te_0$ where $\es_0\oplus\te_0$ is a subalgebra of
$\es\oplus\te$, then
\[
    \der_KL=\ider L\rtimes\de_{\es_0,\te_0}.
\]
\end{thm}
\textsl{Proof:}\newline 1)  It is clear that
$\de_{\m}\subset(\der_K L)_0$.  On the other hand, let
$d=d_{\bz}\oplus d_{\bo}\in(\der_K L)_0$ and so there exist
$M,N\in\ed_K\m,\ \varphi_{\ba}\in\ed_{\bk}\ba,\
\varphi_{\ce}\in\hm_{\bk}(\ba,\ce)$ and
$\psi_{\ba}\in\hm_{\bk}(\ce,\ba),\ \psi_{\ce}\in\ed_{\bk}\ce$
 such that
\[
    d(a,m,n,x)=(\varphi_{\ba}(a)+\psi_{\ba}(x),M(m),N
(n),\varphi_{\ce}(a)+\psi_{\ce}(x)).
\]
Since $d$ is a derivation, for all $a,b\in\ba$ and $x,y\in\ce$ we
have
\begin{align*}
    (\psi_{\ba}([x,y]),&0,0,\psi_{\ce}([x,y]))=d([(a,0,0,x),(b,0,0,y)]\\
    =\ &[(\varphi_{\ba}(a)+\psi_{\ba}(x),0,0,\varphi_{\ce}(a)+\psi_{\ce}(x)),(b,0,0,y)]\\
    &+(-1)^{|(a,0,0,x)||(b,0,0,y)|}[(a,0,0,x),(\varphi_{\ba}(b)+\psi_{\ba}(y),0,
0, \varphi_{\ce}(b)+\psi_{\ce}(y))\\
    =\ &(0,0,0,[\varphi_{\ce}(a)+\psi_{\ce}(x),y]+(-
1)^{|(a,0,0,x)||(b,0,0,y)|}[x,\varphi_{\ce}(b)+\psi_{\ce}(y)])
\end{align*}
which implies that
\[
    \psi_{\ba}[x,y]=0,\ \varphi_{\ce}(\ba)\subset Z(\ce)\textup{ and }
    \psi_{\ce}(\ce)\subset\der_{\bk}\ce.
\]
Moreover, for all $a\in\ba,\ x\in\ce$ and $m\in\m$,
\begin{align*}
    (0,M(am+xm),0,0)&=d(0,am+xm,0,0)=d([(a,0,0,x),(0,m,0,0)])\\
    =\ &[(\varphi_{\ba}(a)+\psi_{\ba}(x),0,0,\varphi_{\ce}(a)+\psi_{\ce}(x),(0,m,0,0)]\\
    &+(-
1)^{|d||a|}[(a,0,0,x),(0,M(m),0,0)]\\
    =\ &(0,(\varphi_{\ba}(a)+\psi_{\ba}(x)+\varphi_{\ce}(a)+\psi_{\ce}(x))(m),0,0)\\
    &+(-
1)^{|d||a+x|}(0,aM(m)+xM(m),0,0).
\end{align*}
Hence
\[
    [M,a\id]=\varphi_{\ba}(a)\id+\varphi_{\ce}(a)\textup{ and }
    [M,x]=\psi_{\ba}(x)\id+\psi_{\ce}(x)
\]
for all $a\in A$ and $x\in\ce$.  Similarly,
\[
    [N,a\id]=-\varphi_{\ba}(a)\id+\varphi_{\ce}(a)\textup{ and }
    [N,x]=-\psi_{\ba}(x)\id+\psi_{\ce}(x)
\]
for all $a\in A$ and $x\in\ce$.  Let $S,T\in\ed_K\m$ be given by
\begin{align*}
    S\ &=\ \half(M+N);\\
    T\ &=\ \half(M-N).
\end{align*}
Then $[S,a\id]=\varphi_{\ba}(a)\id\in A\id$ and
$[T,a\id]=\varphi_{\ce}(a)\in\ce$.  Moreover,
$[S,x]=\psi_{\ce}(x)\in\ce$ and $[T,x]=\psi_{\ba}(x)\id\in\ba\id$.
Hence $S$ satisfies \sones and \sthrees and $T$ satisfies \tones
and \tthree.    Moreover, we have that
\[
    d(a,m,n,x)=([S,a\id]+[T,x],(S+T)m,(S-T)n,[S,x]+[T,a\id])
\]
and so using the fact that $\e_{m,n}=[m^+,n^-]-q(m,n)=-(-
1)^{|m||n|}\e_{n,m}$, we can show that
\begin{eqnarray*}
    (q(S(m),n)+(-
1)^{|m||n|}q(S(n),m))\id&=&[S,q(m,n)\id];\textup{ and}\\
    \e_{T(m),n}+(-
1)^{|m||n|}\e_{T(n),m}&=&[T,q(m,n)\id])
\end{eqnarray*}
and so $S$ and $T$ satisfy \stwos and \ttwos respectively.  Hence
$S\oplus T\in\es\oplus\te$ and so $d\in\de_{\ce}$.  Therefore we
have $\de_{\ce}=(\der_K L)_0$ and so $\de_{\ce}$ is a subalgebra
of $\der_{\bk}\ce_{\infty}$. and so by Proposition \ref{der for
toral L0=A+[L+,L-]}, $\de_{\ce}$ is a subalgebra of $\der_KL$.
Let $\phi: \es \oplus \te \longrightarrow \de_{\ce}$ be defined by
$\phi(S\oplus T)=d_{S,T}$. Clearly this map is well defined and is
\bk-linear since $S$ and $T$ are \bk-linear.  Injectivity follows
from Remark \ref{uniqueness in D} and surjectivity from the
definition of $\de_{\ce}$.  It is easy to check that $\phi$ is a
Lie superalgebra homomorphism and hence isomorphism.\newline 2)
Follows from Proposition \ref{der for toral L0=A+[L+,L-]} and
1).\newline 3)  It is straightforward to verify that
$d_{\ce,\ba\id}\subset\ad L$. Conversely, let $d\in\ider
L\cap\de_{\ce}$.  Then $d=d_{S,T}= \ad(a_0,m_0,n_0,x_0)$ for some
$S\in\es,\ T\in\te,\ (a_0,m_0,n_0,x_0)\in L$.  Then using
(\ref{dformula}) and the fact that $\half\in\bk$, we have that
$S=x_0\in\ce$ and $T=a_0\id\in A\id$.  Hence the assertion
follows.\newline 4)  Assume $\es=\ce\oplus\es_0$ and
$\te=A\id\oplus\te_0$ where $\es_0\oplus\te_0$ is a subalgebra of
$\es\oplus\te$.  Then by Proposition \ref{S+T is LSA},
$\es\oplus\te=(\ce\oplus A\id)\rtimes(\es_0\oplus\te_0).$ Hence by
1), $\{d_{S,T}\mid S\in\es_0,\ T\in\te_0\}$ is a subalgebra of
$\de_{\ce}$ and
\begin{eqnarray*}
    \de_{\ce}&=&\{d_{S,T}\mid S\in\ce,\ T\in A\id\}\rtimes\{d_{S,T}
\mid S\in\es_0,\ T\in\te_0\}\\
    &=&\ad(L_0)\rtimes\{d_{S,T}\mid S\in\es_0,\
T\in\te_0\}.
\end{eqnarray*}
By 3), we have $\ider L\cap\de_{\ce}=\ad(L_0)$. Hence by 2),
\[
    \der_K L=\ider L+\de_{\ce}=\ider L\oplus\{d_{S,T}\mid
S\in\es_0,T\in\te_0\}.
\]
The assertion now follows from 1).\qed\\
\\
In fact, the \ba-supermodules \ess and \tes can be described in
the following terms:  Recall the Jordan superalgebra
$J=\ba\oplus\m$ of Example \ref{JSA of q.f.}.  Given a subalgebra
\ces of \osp(\q) containing \eosp(\q), recall the $J$-supermodule
$X=\ce\oplus\m$ with action of $J$ on $X$ given by (\ref{jsm}).
Let
\begin{eqnarray*}
    \der_*(J)&=&\{d\in\der_{\bk}(J)\mid
d(\ba)\subset\ba\textup{ and }d(\m)\subset\m\};\textup{ and}\\
    \der_*(J,X)&=&\{d\in\der_{\bk}(J,X)\mid d(\ba)\subset\ce\textup{
and }d(\m)\subset\m\}.
\end{eqnarray*}
\begin{prop}\label{S and T as derivations}
In the above setting we have that
\begin{itemize}
\item[(1)]  the map $\der_*(J)\rightarrow\esone\cap\estwo:d\mapsto
d|_{\m}$ is a Lie superalgebra isomorphism; and \item[(2)]  the
map $\der_*(J,X)\rightarrow\teone\cap\tetwo:d\mapsto d|_{\m}$ is a
\bk-supermodule isomorphism.
\end{itemize}
\end{prop}
\textsl{Proof:}  (1)  Let $d\in\ed_K\m$ such that
$d(\ba)\subset\ba$ and $d(\m)\subset\m$.  Then from the definition
of $\der(J)$ and Proposition \ref{prop of S and T}.1.(i) it
follows that $d\in\der_*(J)$ if and only if $d|_{\m}$ satisfies
\sones and \stwo.\newline (2)  Let
$d\in\hm_{\bk}(\ba\oplus\m,\ce\oplus\m)$ such that
$d(\ba)\subset\ce$ and $d(\m)\subset\m$.  Then from the definition
of $\der(J,X)$ and Proposition \ref{prop of S and T}.2.(i) it
follows that $d\in\der_*(J,X)$ if and only if $d|_{\m}$ satisfies
\tones and \ttwo.\qed


\section{DERIVATIONS AS A SEMIDIRECT PRO\-DUCT}
Here we will present some examples of subalgebras of \osp(\q)
containing \eosp(\q) where we get a direct splitting of the
algebra of derivations into inner and outer
derivations.  We keep the setting of Section 3.\\
\\
We let \bastars to be the set of all invertible homogeneous
elements of \ba.  Note that, since $\half\in\ba$, $\bastar\subset\ba_{\bz}$.\\
\\
Given a free \ba-supermodule $\m=\oplus_{i\in I}\ba m_i$ with a
quadratic form \q, we say that \qs is \emph{almost diagonalizable}
if for each $i\in I$ there exists $\und{i}\in I$ such that
\begin{align*}
    &q(m_i,m_{\und{i}})\in\bastar;\textup{ and}\\
    &q(m_i,m_j)=0\textup{ for all } j\in I,\,j\neq\und{i}.
\end{align*}
If \qs is almost diagonalizable then for each $i\in I$, $\und{i}$
is unique, $\und{\und{i}}=i$ and $|m_i|=|m_{\und{i}}|$.  In
addition, if $q$ is \emph{almost diagonalizable} with respect to
$\{m_i\mid i\in I\}$ where $i=\und{i}$ then we say that \qs is
\emph{invertibly diagonalizable}.  Note that if \qs is invertibly
diagonalizable with respect to $\{m_i\mid i\in I\}$ then the basis
elements $m_i$ are even.
\begin{prop}\label{T for even and q=inv diag+0}\
Let $\p$ be an \ba-supermodule with an \ba-quadratic form $q_{\p}$
and an element $p\in\p$ such that $\ann_{\ba}(p)=(0)$.  Let
$\n=\oplus_{i=\pm 1}\ba n_i$ be a free \ba-supermodule with an
invertibly diagonalizable \ba-quadratic form $\q_{\n}$. Set
$\m=\n\oplus\p$ with $q=q_{\n}\oplus q_{\p}$ and let \ces be a
subalgebra of \osp(\q) containing \eosp(\q).  Then
\[
    \te_{\ce}=\tetwo=\ba\id.
\]
\end{prop}
\textsl{Proof:} Recall that $\ba\id\subset\te\subset\tetwo$.  Let
$T\in\tetwo$ and write
\[
    T=\left(\begin{array}{cc}
    T_{\n,\n}&T_{\p,\n}\\
    T_{\n,\p}&T_{\p,\p}
    \end{array}\right)
\]
with respect to $\left(\begin{array}{c}\n\\ \p\end{array}\right)$,
so $T_{\n,\p}\in\hm_{\bk}(\n,\p)$, etc. Let $p\in\p$ and let
$i,k\in \{\pm 1\},\,k\neq{i}$.  Then by \ttwo,
\begin{eqnarray}
\nonumber   0&=&[T,q(an_i,p)\id](n_k)=
\e_{T(an_i),p}(n_k)+(-1)^{|a||p|}\e_{T(p),an_i}(n_k)\\
\label{T(N,N)}  &=&-q(T(an_i),n_k)p-(-1)^{|a||p|}q(T(p),n_k)an_i.
\end{eqnarray}
Since $\m=\n\oplus\p$ and $q(\n,\p)=q(\p,\n)=0$ we have that if
$p\in\p$ is such that $\ann_{\ba}(p)=(0)$, then
\[
    0=q(T(an_i),n_k)=q_{\n}(T_{\n,\n}(an_i),n_k) \textup{ for all
}i,k\in I,\,k\neq{i}.
\]
We write $T_{\n,\n}(an_i)=t_{i,1}(a)n_1+t_{i,-1}(a)n_{- 1}$ where
$t_i\in\ed_{\bk}\ba$ for $i\in \{\pm 1\}$ and we obtain
\[
    0=t_{i,{k}}(a)q(n_{{k}},n_k).
\]
Since $q(n_k,n_k)\in\bastar$, we therefore have that
$t_{i,k}(a)=0$ for $i\neq k$ and so
$T_{\n,\n}(an_i)=t_{ii}(a)n_i$.  Moreover, since $p\in\p$ and
$i\in\{\pm 1\}$ were arbitrary in (\ref{T(N,N)}) we have that
\[
    0=q(T(p),n)=q_{\n}(T_{\p,\n}(p),n)\textup{ for
all }n\in\n,\,p\in\p
\] and so since $q_{\n}$ is nondegenerate, we have that $T_{\p,\n}=0$.
Now, let $i\in \{\pm 1\}$. Then for all $p\in\p$,
\begin{eqnarray*}
    0&=&[T,q(n_i,p)\id](n_{{i}})=\e_{T(n_i),p}(n_{{i}})+
\e_{T(p),n_i}(n_{{i}})\\
    &=&-t_{ii}(1)q(n_i,n_{{i}})p-T(p)
q(n_i,n_{{i}})
\end{eqnarray*}
and so we have that $T(p)=t_{ii}(1)p,\ i\in\{\pm 1\}$.  Set
$t:=t_{ii}(1)$.\newline\newline Let $i,k\in \{\pm 1\},\,i\neq {k}$
and let $a\in\ba$.  Then
\begin{eqnarray*}
    0&=&[T,q(n_{{i}},an_k)\id](n_i)=(\e_{T(n_{{i}}),an_k}+
\e_{T(an_k),n_{{i}}})(n_i)\\
    &=&-tq(n_{{i}},n_i)an_k-t_{kk}(a)n_kq(n_{{i}},n_i)-T_{\n,\p}(an_k)q(n_{{i}},n_i).
\end{eqnarray*}
Since $q(n_{{i}},n_i)\in\bastar$, $T_{\n,\p}=0$ and
$t_{kk}(a)=ta$.  Hence $T=t\id$ and so $\te=\ba\id$.\qed
\begin{rmk}\textup{  For $I$ and \ps different from the above the
situation seems to be more complicated.}
\end{rmk}
One can easily show that if \ns and \ps are $A$-supermodules with
\ba-quadratic forms $q_{\n}$ and $q_{\p}$ respectively, then for
$\m=\n\oplus\p$ with $q=q_{\n}\oplus q_{\p}$,
\begin{eqnarray*}
\eosp(q)=\eosp(q_{\n})\oplus\eosp(q_{\n},q_{\p})\oplus\eosp(q_{\p})
\end{eqnarray*}
where $\eosp(q_{\n},q_{\p})=\textup{span}_A\{\e_{n,p}\mid n\in\n,\
p\in\p\}$ and
\begin{eqnarray*}
\ovl{\eosp(q)}=\ovl{\eosp(q_{\n})}\oplus\ovl{\eosp(q_{\n},q_{\p})}
\oplus\ovl{\eosp(q_{\p})}.
\end{eqnarray*}
Moreover, if $q_{\n}$ is invertibly diagonalizable then
\begin{eqnarray*}
\osp(q)=\ovl{\eosp(q_{\n})}\oplus\ovl{\eosp(q_{\n},q_{\p})}\oplus\osp(q_{\p}).
\end{eqnarray*}
Given $S\in\esone$ we denote by $\ad_{\ba}S$ the derivation
$\Delta\in\der_{\bk}\ba$ which satisfies $[S,a\id]=\Delta(a)\id$
for all $a\in\ba$.  We define
\[
    \es_{\n}^{\p}=\{S_{\n}\oplus
S_{\p}\in(\esone_{\n}\cap\estwo_{\n})\oplus(\esone_{\p}\cap\estwo_{\p})\mid
\ad_{\ba}S_{\n}=\ad_{\ba}S_{\p}\}.
\]
\begin{prop}\label{S for N+P}
Let $\m=\n\oplus\p$ where \ns and \ps are \ba-supermodules with
\ba-quadratic forms $\q_{\n}$ and $q_{\p}$ respectively.  Assume
that $q_{\n}$ is almost diagonalizable with respect to some basis
$\{n_i\mid i\in I\}$.  Then
\begin{eqnarray*}
    \esone\cap\estwo&=&\ovl{\eosp(q_{\n},q_{\p})}\oplus\es_{\n}^{\p}
\end{eqnarray*}
and
\begin{eqnarray*}
    \es&=&\{S\in\ovl{\eosp(q_{\n},q_{\p})}\oplus\es_{\n}^{\p}\mid
[S,\ce]\subset\ce\}.
\end{eqnarray*}
\end{prop}
\textsl{Proof:}  Let $S\in\esone\cap\estwo$. By \sone,
$[S,a\id]=\Delta(a)\id$ for some $\Delta\in\der_{\bk}\ba$. Write
\[
    S=\left(\begin{array}{cc}S_{\n,\n}&S_{\n,\p}\\S_{\p,\n}&S_{\p,\p}
\end{array}\right)
\]
where $S_{\n,\p}\in\hm_{\bk}(\n,\p)$, etc.  Then using \sone,
$[S_{\n,\n}\oplus S_{\p,\p},a\id]=\Delta(a)\id$ and
$S_{\n,\p}\oplus S_{\p,\n}\in\hm_{\ba}\m.$  In addition, using the
orthogonality of \ns and \ps and \sones applied to \ns and \ps
separately we have $S_{\n,\n}\oplus S_{\p,\p}\in\es_{\n}^{\p}$.
Moreover, using \stwos we can show that
\[
    q((S_{\n,\p}\oplus S_{\p,\n})(m),m')+(-
1)^{|m||m'|}q((S_{\n,\p}\oplus S_{\p,\n})(m'),m)=0
\]
and so $S_{\n,\p}\oplus
S_{\p,\n}\in\osp(q)=\ovl{\eosp(q_{\n})}\oplus
\ovl{\eosp(q_{\n},q_{\p})}\oplus\osp(q_{\p})$. Hence clearly
$S_{\n,\p}\oplus S_{\p,\n}\in\ovl{\eosp(q_{\n},q_{\p})}$ and so
$\es^{(1)}\cap\es^{(2)}\subset\ovl{\eosp(q_{\n},q_{\p})}\oplus
\es^{\p}_{\n}.$  Conversely, to show that
$\ovl{\eosp(q_{\n},q_{\p})}\oplus\es^{\p}_{\n}\subset\esone\cap\estwo$,
by Proposition \ref{prop of S and T}.1.4, we only need to show
that $\es^{\p}_{\n}\subset\esone\cap\estwo$.  Let $S_{\n}\oplus
S_{\p}\in\es_{\n}^{\p}$ and let $\Delta\in\der_{\bk}\ba$ such that
$[S_{\n}\oplus S_{\p},a\id]=\Delta(a)\id$.  Hence clearly
$S_{\n}\oplus S_{\p}\in\esone$.  One can easily show that
$S_{\n}\oplus S_{\p}\in\estwo$ and so the rest of the assertion
follows.\qed
\begin{prop}\label{S for s-orthogonal}
Let \ns and \calrs be free with homogeneous bases $\{n_i\mid i\in
I\}$ and $\{r_j\mid j\in J\}$ and \ba-quadratic forms $\q_{\n}$
and $q_{\calr}$ respectively where $q_{\n}$ is almost
diagonalizable with respect to $\{n_i\mid i\in I\}$ and
$q_{\calr}\equiv 0$. Set $\m=\n\oplus\calr$ with the \ba-quadratic
form $q=q_{\n}\oplus q_{\calr}$. For $\Delta\in\der_{\bk}\ba$
define $\Delta_{\m}\in\ed_{\ba}\m$ by setting
\begin{align}
\nonumber   \Delta_{\m}\left(\sum_{i\in I}a_in_i+\sum_{j\in
J}b_jr_j\right)=\ &\sum_{i\in
I}\left(\half\Delta(q(n_i,n_{\und{i}}))q(n_i,n_{\und{i}})\inv
a_i+\Delta(a_i)\right)n_i\\
\label{delta_M} &+\sum_{j\in J}\Delta(b_j)r_j.
\end{align}
Set
$(\der_{\bk}\ba)_{\m}=\{\Delta_{\m}\mid\Delta\in\der_{\bk}\ba\}$.
Then the map $\der_{\bk}\ba\rightarrow
(\der_{\bk}\ba)_{\m}:\Delta\mapsto\Delta_{\m}$ is a Lie
superalgebra isomorphism and
\begin{eqnarray*}
    \esone\cap\estwo&=&\osp(q)\rtimes(\der_{\bk}\ba)_{\m}\\
    &=&\ovl{\eosp(q)}\rtimes(\ed_{\ba}\calr\rtimes(\der_{\bk}\ba)_{\m}).
\end{eqnarray*}
Hence
\[
    \es_{\ce}=\{S\in\osp(q)\rtimes(\der_{\bk}\ba)_{\m}
\mid [S,\ce]\subset\ce\}
\]
and in particular, if $\ce=\eosp(\q),\ \ovl{\eosp(q)}$ or \osp(\q)
then
\begin{eqnarray*}
    \es_{\ce}&=&\esone\cap\estwo=\osp(q)\rtimes (\der_{\bk}\ba)_{\m}\\
    &=&\ovl{\eosp(q)}\rtimes(\ed_{\ba}\calr\rtimes(\der_{\bk}\ba)_{\m}).
\end{eqnarray*}
\end{prop}
\textsl{Proof:}  Let $S\in\esone\cap\estwo$. By \sone,
$[S,a\id]=\Delta(a)\id$ for some $\Delta\in\der_{\bk}\ba$. Then
for each $m\in\m$ and $a\in \ba$, $S(am)=(-
1)^{|S||a|}aS(m)+\Delta(a)m.$  Set $S_0\in\ed_{\ba}\m$ by
\begin{eqnarray*}
    S_0(n_i)=S(n_i)-
\half\Delta(q(n_i,n_{\und{i}}))q(n_i,n_{\und{i}})\inv
n_i\textup{\quad and\quad }
    S_0(r_j)=S(r_j)
\end{eqnarray*}
for $i\in I,\ j\in J$.  Then for each $i,k\in I$, using \stwo,
$q(S_0(n_i),n_k)+(- 1)^{|n_i||n_k|}q(S_0(n_k),n_i)=0.$ Moreover,
for each $i\in I,\, j\in J$,
\begin{align*}
    q(S_0(r_j),n_i)+(-
1)^{|r_j||n_i|}q(S_0(n_i),r_j)=q(S(r_j),n_i)
    \overset{\stwo}{=}[S,q(r_j,n_i)\id]=0
\end{align*}
and for all $j,l\in J$,
\[
    q(S_0(r_j),r_l)+(-
1)^{|r_j||r_l|}q(S_0(r_l),r_j)=0.
\]
Hence since $S_0\in\ed_{\ba}\m$, $S_0\in\osp(q)$ and so by a
direct calculation we get that $S=S_0+\Delta_{\m}$.  Therefore
$\esone\cap\estwo\subset\osp(q)+(\der_{\bk}\ba)_{\m}$. Conversely,
to show that
$\osp(q)+(\der_{\bk}\ba)_{\m}\subset\esone\cap\estwo$, by
Proposition \ref{prop of S and T}.1.4, we only need to show that
$(\der_{\bk}\ba)_{\m}\subset\esone\cap\estwo$.  This again follows
from a direct calculation.  Since
$(\der_{\bk}\ba)_{\m}\cap\ed_{\ba}\m=(0)$ and
$\osp(q)\in\ed_{\ba}\m$,
$\esone\cap\estwo=\osp(q)\oplus(\der_{\bk}\ba)_{\m}$.  One can
check that the map
$\der_KA\rightarrow(\der_KA)_{\m}:\Delta\mapsto\Delta_{\m}$ is a
Lie superalgebra isomorphism. Hence $(\der_KA)_{\m}$ is a
subalgebra of $\esone\cap\estwo$ and so since by Proposition
\ref{prop of S and T}.1.(iv) \osp(\q) is an ideal of
$\esone\cap\estwo$, it follows that
$\esone\cap\estwo=\osp(q)\rtimes(\der_{\bk}\ba)_{\m}.$ Since
$q_{\calr}=0$, $\osp(q)=\ovl{\eosp(q)}\rtimes\ed_{\bk}\calr$ and
so
\[
    \esone\cap\estwo=(\ovl{\eosp(q)}\rtimes\ed_{\bk}\calr)\rtimes(\der_{\bk}\ba)_{\m}=
\ovl{\eosp(q)}\oplus(\ed_{\bk}\calr\oplus(\der_{\bk}\ba)_{\m})
\]
Since $(\der_{\bk}\ba)_{\m}$ and $\ed_{\ba}\calr$ are subalgebras
of $\esone\cap\estwo$ to show that
$\ed_A\calr\oplus(\der_KA)_{\m}$ is a subalgebra of
$\esone\cap\estwo$, we only need to show that
$[(\der_KA)_{\m},\ed_A\calr]\subset
(\der_KA)_{\m}\oplus\ed_A\calr$.  In fact, we will show that
\[
    [(\der_KA)_{\m},\ed_A\calr]\subset\ed_A\calr.
\]
Indeed, let $\Delta\in\der_KA$ and $S\in\ed_A\calr$.  Then with
respect to $\m=\left(\begin{array}{c}\n\\
\calr\end{array}\right)$, we have
\[
    [\Delta_{\m},S]=\left[\left(
\begin{array}{cc}
    \Delta|_{\n}&0\\
    0&\Delta|_{\calr}
\end{array}\right),\left(
\begin{array}{cc}
    0&0\\
    0&S
\end{array}\right)\right]=\left(
\begin{array}{cc}
    0&0\\
    0&[\Delta|_{\calr},S]
\end{array}\right)\in\hm_{\bk}(\m,\calr)
\]
and so $[(\der_KA)_{\m},\ed_A\calr]\subset\hm_{\bk}(\m,\calr)$.
Moreover, if $S(r_j)=\sum_{i\in J}s_{ji}r_i$ for some
$s_{ji}\in\ba,\,i,j\in J$ then $[[\Delta_{\m},S],a\id](br_j)=0$
and so $[(\der_KA)_{\m},\ed_A\calr]\subset
0\oplus\ed_A\calr=\ed_A\calr.$ Hence
$[(\der_KA)_{\m}\oplus\ed_A\calr,(\der_KA)_{\m}\oplus\ed_A\calr]
\subset(\der_KA)_{\m}\oplus\ed_A\calr$ and so
$(\der_KA)_{\m}\oplus\ed_A\calr$ is a subalgebra of $\es$. The
rest of the assertion now follows from the definition of \ess and
from
Proposition \ref{prop of S and T}.1.(iii).\qed\\
\\
For what follows we need the following:
\begin{Defn}\textup{
We define the \emph{locally inner derivations of L} to be the
derivations $d\in\der_KL$ such that for any finite subset
$\{x_1,\ldots, x_n\}$ of $L$ there exists $x\in L$ satisfying
$d(x_i)=[x,x_i]$ (see \cite{DZ}). We denote the locally inner
derivations of $L$ by $\adloc L$. }\end{Defn}
\begin{rmk}\textup{
One can easily show that $\adloc L$ is an ideal of $\der_{\bk} L$
and moreover that $\ovl{\eosp(q)}\subset\adloc\eosp(q_{\infty})$
(one uses the fact that
$[x,\e_{m,n}]=\e_{x(m),n}+(-1)^{|x||m|}\e_{m,x(n)}$ for
$x\in\osp(q)$ and $m,n\in\m$). }\end{rmk} We will, given
$S\in\esone$, for simplicity write $[S,\ba\id]$ for the map
$\Delta\in\der_{\bk}\ba$ satisfying $[S,a\id]=\Delta(a)\id$.
\begin{co}\label{S for almost diagonalizable+P}
Let \ns and \ps be \ba-supermodules with \ba-quadratic forms
$q_{\n}$ and $q_{\p}$ respectively where $q_{\n}$ is almost
diagonalizable with respect to some basis $\{n_i\mid i\in I\}$ of
\n.  Then
\begin{align*}
\begin{array}{rcl}
    \ovl{\eosp(q_{\n})}\oplus\ovl{\eosp(q_{\n},q_{\p})}\oplus(\esone_{\p}\cap\estwo_{\p})
&\rightarrow&\esone\cap\estwo\\
    x\oplus S&\mapsto& x\oplus (\ad_{\ba} S)_{\n}\oplus S
\end{array}
\end{align*}
is a Lie superalgebra isomorphism. If $\ce=\eosp(q),\
\ovl{\eosp(q)}$ or \osp(\q) then
\[
    \es\cong \ovl{\eosp(q_{\n})}\oplus\ovl{\eosp(q_{\n},q_{\p})}\oplus\es_{\p}
\]
and in particular, if \qs is the orthogonal sum of an almost
diagonalizable quadratic form and an invertibly diagonalizable
\ba-quadratic form of rank 2, then
\[
    \der_{\bk}\ce_{\infty}\cong\left\{\begin{array}{ll}\adloc(\ce_{\infty})+\es_{\p},& i\!f\
\ce=\eosp(q);\\
        \ad(\ce_{\infty})+\es_{\p},&i\!f\ \ce=\ovl{\eosp(q)},\
\osp(q).\end{array}\right.
\]
\end{co}
\textsl{Proof:}  By Proposition \ref{S for N+P},
$\esone\cap\estwo=\ovl{\eosp(q_{\n},q_{\p})}\oplus \es_{\n}^{\p}.$
By Proposition \ref{S for s-orthogonal},
$\es_{\n}=\ovl{\eosp(q_{\n})}\rtimes(\der_{\bk}\ba)_{\n}$.  We
claim that
\begin{eqnarray}\label{splitting of S(N,P)}
    \es_{\n}^{\p}=(\ovl{\eosp(q_{\n})}\oplus
0)\oplus\{(\ad_{\ba}S_{\p})_{\n}\oplus S_{\p}\mid
S_{\p}\in\es_{\p}\}.
\end{eqnarray}
Indeed, if $S_{\n}\oplus S_{\p}\in\es_{\n}^{\p}$ let
$\Delta=\ad_{\ba}S_{\n}=\ad_{\ba}S_{\p}$.  Then
$S_{\n}=x+\Delta_{\n}\in\ovl{\eosp(q_{\n})}\oplus(\der_{\bk}\ba)_{\n}$
and so $S_{\n}\oplus S_{\p}=(x\oplus 0)+(\Delta_{\n}\oplus
S_{\p})=(x\oplus 0)+(([S_{\p},\ba\id])_{\n}\oplus S_{\p})$ which
is the element of the right-hand side of (\ref{splitting of
S(N,P)}). Conversely, since
$\ovl{\eosp(q_{\n})}\subset\ed_{\ba}\m$ we have
$\ovl{\eosp(q_{\n})}\subset\es_{\n}$ with $\ad_{\ba}S_{\n}=0$.
Moreover, given $S_{\p}\in\es_{\p}$ with
$\ad_{\ba}S_{\p}=\Delta\in\der_{\bk}\ba$ we have
$\Delta_{\n}\oplus S_{\p}\in\es_{\n}^{\p}$ since $\Delta_{\n}\in
S_{\n}$ by Proposition \ref{S for s-orthogonal}.  If $0=(x\oplus
0)+((\ad_{\ba}S_{\p})_{\n}\oplus S_{\p})$ then $S_{\p}=0$ hence
$(\ad_{\ba}S_{\p})_{\n}=0$ and so $x=0$.  Therefore
\[
    \esone\cap\estwo=\ovl{\eosp(q_{\n},q_{\p})}\oplus\ovl{\eosp(q_{\n})}
\oplus\{(\ad_{\ba}S_{\p})_{\n}\oplus S_{\p}\mid
S_{\p}\in\es_{\p}\}.
\]
It is easy to see that $\es_{\p}\rightarrow
\{(\ad_{\ba}S_{\p})_{\n}\oplus S_{\p}\mid
S_{\p}\in\es_{\p}\}:S_{\p}\mapsto (\ad_{\ba}S_{\p})_{\n}\oplus
S_{\p}$ is a Lie superalgebra isomorphism and so the first
assertion of the Corollary holds.  The second assertion follows
from \ref{prop of S and T}.1.(iii).  If $\ce=\eosp(q),\
\ovl{\eosp(q})$ or $\osp(q)$ then by Proposition \ref{T for even
and q=inv diag+0}, $\es=\esone\cap\estwo$.  The rest of the
Corollary now follows from Theorem \ref{main} using the fact that
$\ovl{\eosp(q_{\n},q_{\p})}\oplus\ovl{\eosp(q_{\n})}\subset\ovl{\eosp(q)}$
and that $\ovl{\eosp(q)}\subset \adloc(\eosp(q_{\infty}))$.\qed
\begin{rmk}\textup{
In \cite{Be}, Benkart described the derivations of root-graded Lie
algebras over a field \F\  of characteristic $0$.  In particular,
for a centreless Lie algebra $L$ graded by finite root systems $B$
or $D$ she proves that
\[
    \der_{\F}L=\ider L+\der_*(J)
\]
for a certain Jordan algebra $J$.  Since $L$ is isomorphic to
$\ce=\eosp(\q_{I}\oplus q_{x_0}\oplus q_{\p})$ for some finite
index set $I$, base point $x_0$ and an \ba-module \ps where \bas
is an extension of \F\ and $\q_{I}\oplus q_{x_0}$ is almost
diagonalizable.  In addition, since $I$ is finite, one can show
that $\ovl{\eosp(q_{\n},q_{\p})}\oplus
\ovl{\eosp(q_{\n})}=\eosp(q_{\n},q_{\p})\oplus \eosp(q_{\n})$ and
so using the above corollary,
$\der_{\bk}\ce_{\infty}\cong\ad(\ce_{\infty})+\es_{\p}$.  However,
by Proposition \ref{S and T as derivations},
$\es_{\p}\cong\der_*(A\oplus\p)$ where $A\oplus\p$ is a Jordan
algebra and this algebra is isomorphic to the Jordan algebra which
appears in Benkart's description of $B$-graded Lie algebras.
}\end{rmk}
\begin{thm}\label{der for q=inv diag+0}
Let $\n$ and $\calr$ be free \ba-supermodules with
$\textup{rank}(\n\oplus\calr)>2$ and \ba-quadratic forms $q_{\n}$
and $q_{\calr}$ respectively such that $q_{\calr}\equiv 0$ and
$q_{\n}$ is an orthogonal sum of an almost diagonalizable
\ba-quadratic form and an invertibly diagonalizable \ba-quadratic
form on a free \ba-supermodule of rank 2.  Let $\m=\n\oplus\calr$
with the \ba-quadratic form $q=q_{\n}\oplus q_{\calr}$ and let
$\ce$ be a subalgebra of \osp(\q) containing \eosp(\q).  For
$\Delta\in\der_{\bk}\ba$ define $\Delta_{\m}\in\ed_{\bk}\m$ using
(\ref{delta_M}). Let the map $\{S\in\osp(q)\oplus
(\der_KA)_{\m}\mid
[S,\ce]\subset\ce\}\rightarrow\de_{\ce}:S\mapsto d_{S,0}$ be the
Lie superalgebra monomorphism given by (\ref{dformula}).  Then
\(\der_KA\rightarrow(\der_KA)_{\m}:\Delta\mapsto \Delta_{\m}\) is
a Lie superalgebra isomorphism and
\[
    \der_{\bk}\ce_{\infty}=\ad\ce_{\infty}+
\{d_{S,0}\in\de_{\ce}\mid S\in\osp(q)\oplus(\der_KA)_{\m}\text{
such that }[S,\ce]\subset\ce\}.
\]
with
\begin{eqnarray*}
    \ad\ce_{\infty}\cap\{d_{S,0}\in\de_{\ce}\mid S\in\osp(q)\text{ such that
}[S,\ce]\subset\ce\}=\ad\ce=\{d_{S,0}\mid S\in\ce\}.
\end{eqnarray*}
Hence, whenever $\es=\ce\rtimes\es_0$,
\[
    \der_{\bk}\ce_{\infty}=\ad(\ce_{\infty})\rtimes\{d_{S,0}\in\de_{\ce}\mid
S\in \es_0\}.
\]
In particular:
\begin{itemize}
\item[(1)]  If $\ce=\eosp(\q)$ then
\[
    \der_{\bk}\ce_{\infty}=\adloc\ce_{\infty}\rtimes
\{d_{S,0}\in\de\mid
S\in\ed_{\ba}\calr\oplus(\der_{\bk}\ba)_{\m}\}.
\]
Moreover, if \ms is free of finite rank then
\[
    \der_{\bk}\ce_{\infty}=\ad(\ce_{\infty})\rtimes
\{d_{S,0}\in\de\mid
S\in\ed_{\ba}\calr\oplus(\der_{\bk}\ba)_{\m}\}.
\]
\item[(2)]  If $\ce=\ovl{\eosp(q)}$ then
\[
    \der_{\bk}\ce_{\infty}=\ad(\ce_{\infty})\rtimes
\{d_{S,0}\in\de\mid
S\in\ed_{\ba}\calr\oplus(\der_{\bk}\ba)_{\m}\}.
\]
\item[(3)]  If $\ce=\osp(q)$ then
\[
    \der_{\bk}\ce_{\infty}=\ad\ce_{\infty}\rtimes
\{d_{S,0}\in\de\mid S\in(\der_{\bk}\ba)_{\m}\}.
\]
\end{itemize}
\end{thm}
\textsl{Proof:}  By Theorem \ref{main}.2 we have that the map
$\es\oplus\te\rightarrow\de_{\ce}:(S,T)\mapsto d_{S,T}$ is a Lie
superalgebra isomorphism with the induced bracket operation and
$\der_{\bk}\ce_{\infty}=\ad(\ce_{\infty})+\de_{\ce}.$ By
Proposition \ref{S for s-orthogonal},
$\es=\{S\in\osp(q)\rtimes(\der_KA)_{\m}\mid [S,\ce]\subset\ce\}$
where the map
$\der_KA\rightarrow(\der_KA)_{\m}:\Delta\mapsto\Delta_{\m}$ is a
Lie superalgebra isomorphism.  Moreover, by Proposition \ref{T for
even and q=inv diag+0}, $\te=\ba\id$.  The main assertion now
follows from Theorems \ref{main}.3 and \ref{main}.4. and so (3)
follows immediately from this and (2) follows from the fact that
$\osp(\q)=\ovl{\eosp(q)}\rtimes \ed_{\ba}\calr$ in this setting.
Let $\ce=\eosp(q)$.  By Proposition \ref{S for s-orthogonal},
\begin{eqnarray*}
    \es=\osp(q)\rtimes(\der_{\bk}\ba)_{\m}=\ovl{\eosp(\q)}\rtimes(\ed_{\ba}\calr\oplus
(\der_{\bk}\ba)_{\m})
\end{eqnarray*}
and so
\[
    \der_{\bk}\ce_{\infty}=\ad\ce_{\infty}+
(\ovl{\eosp(q)}\oplus\{d_{S,0}\in\de_{\ce}\mid
S\in\ed_{\ba}\calr\oplus(\der_{\bk}\ba)_{\m}\}.
\]
We will show that
$\adloc(\ce_{\infty})=\ad(\ce_{\infty})+\ovl{\eosp(q)}$. The
inclusion $\supset$ is clear.  Since
$\adloc(\ce_{\infty})\subset\ed_{\ba}\ce_{\infty}$, from
(\ref{dformula}) we have $\adloc(\ce_{\infty})\cap
d_{\ed_{\ba}\calr\oplus(\der_{\bk}
\ba)_{\m},0}=\adloc(\ce_{\infty})\cap d_{\ed_{\ba}\calr,0}$.
Suppose that $d_{S,0}\in\adloc(\eosp(q_{\infty}))$ for some
$S\in\ed_{\ba}\calr$.  Then for each $j\in R$ there exists
$(a_0,m_0,n_0,x_0)\in\eosp(q_{\infty})$ such that
\begin{eqnarray*}
    d_{S,0}(0,r_j,0,0)&=&(0,S(r_j),0,0)=[(a_0,m_0,n_0,x_0),
(0,r_j,0,0)]=(0,a_0r_j,0,\e_{n_0,r_j})\\
    d_{S,0}(0,0,r_j,0)&=&(0,0,S(r_j),0)=[(a_0,m_0,n_0,x_0),
(0,0,r_j,0)]=(0,0,-a_0r_j,\e_{m_0,r_j})
\end{eqnarray*}
and so $S(r_j)=a_0r_j$ and $S(r_j)=-a_0r_j$ hence $S(r_j)=0$.
Since $j\in R$ was arbitrary, $S=0$ and so
$\adloc(\ce_{\infty})\cap d_{\ed_{\ba}\calr,0}=(0)$. Hence since
clearly $\ad\ce_{\infty}\subset\adloc\ce_{\infty}$ and since
$\adloc (\ce_{\infty})$ is an ideal of $\ce_{\infty}$,
\[
    \der_K(\ce_{\infty})=\adloc(\ce_{\infty})\rtimes
d_{\ed_A\calr\oplus(\der_KA)_{\m},0}.
\]
Now, if \ms is of finite rank then one can show that
$\ovl{\eosp(q)}=\eosp(q)$ (see \cite{D}) and so
$\adloc\ce_{\infty}=\ad\ce_{\infty}$
hence the assertion follows.\qed\\
\\
The setting of Theorem \ref{der for q=inv diag+0} also holds in
the following example:  Let $\m=\oplus_{i\in I}\ba m_i$ be a free
supermodule where $I=\Z$ or $I=\half\Z$ and let \qs be the
quadratic form on \ms defined by
\begin{eqnarray*}
    q(m_i,m_j)&=&(-1)^i\delta_{i,-j},\quad i,j\in\Z;\\
    q(m_i,m_j)&=&(-1)^{i+\half}\delta_{i,-j},\quad i,j\in\half+\Z;\\
    q(m_i,m_j)&=&0,\quad i\in\Z,\,j\in\half+\Z.
\end{eqnarray*}
Define $\ospfd(q)$ to be the set of all \ba-endomorphisms
$(a_{ij})$ in \osp(\q) of \ms which have finitely many nonzero
diagonals.  For $\bk=\ba=\C$, the algebra $\ospfd(q)$ has been
studied by Kac (\cite{Kac1}) for $I=\Z$ and by Cheng-Wang
(\cite{CW}) for $I=\half\Z$.
\begin{co}\label{der of Kac}
If the elements of $\Z$ are invertible in \bas then \ospfd(\q) is
self-normalizing in \osp(\q) and
\[
    \der_{\bk}\ce_{\infty}=\ad(\ce_{\infty})\rtimes\{d_{S,0}\in\de_{\ce}\mid
S\in(\der_{\bk}\ba)_{\m}\}.
\]
\end{co}
\textsl{Proof:}  By Proposition \ref{S for s-orthogonal},
\[
    \es=\{S\in\osp(q)\oplus(\der_{\bk}\ba)_{\m}\mid[S,\ce]\subset\ce\}.
\]
Since the elements of $(\der_{\bk}\ba)_{\m}$ satisfy \sone, by
Proposition \ref{prop of S and T}.1.1,
$[(\der_{\bk}\ba)_{\m},\allowbreak\ed_{\ba}\m]\subset\ed_{\ba}\m$.
Hence since the elements of $(\der_{\bk}\ba)_{\m}$ are diagonal,
$[(\der_{\bk}\ba)_{\m},\allowbreak\ce]\subset\ce$.  Let $\{v_i\mid
i\in I\cup J\}$ denote the standard basis of \m.  Let
$(a_{ij})_{i,j\in I\cup J}\in\osp(q)$.  Since \ms is almost
diagonalizable, one can show that $a_{i,j}=-a_{-j,-i}q(v_i,v_{-
i})q(v_j,v_{-j})$ for all $i,j\in I\cup J$.  Define
$(b_{ij})_{i,j\in I\cup J}\in\ospfd(q)$ by $b_{k,k}=k,\ k\neq 0$
and $b_{0,0}=0$.  Then $[(a_{ij}),(b_{ij})]\in\ospfd{q}$ implies
that $(a_{ij})\in\ospfd(q)$ and so $\ospfd(q)$ is
self-normalizing.  Hence
\begin{eqnarray*}
    \es&=&\{S\in\osp(q)\rtimes(\der_{\bk}\ba)_{\m}\mid[S,\ce]\subset\ce\}=\{S\in\osp(q)
\mid[S,\ce]\subset\ce\}\rtimes(\der_{\bk}\ba)_{\m}\\
    &=&\ospfd(q)\rtimes(\der_{\bk}\ba)_{\m}
\end{eqnarray*}
and so
$\der_{\bk}\ce_{\infty}=\ad(\ce_{\infty})+\{d_{S,0}\in\de_{\ce}\mid
S\in\ospfd(q)\rtimes(\der_{\bk}\ba)_{\m}\}.$
The assertion now follows from Theorem \ref{main}.4.\qed\\
\\
The following corollary is a generalization of Benkart's result
for Lie algebras of type $B$ and $D$ (\cite[Theorem 3.6]{Be}).
\begin{co}
Suppose $\m=H(I;\ba)\oplus\p$ where $|I|\geq 2$ or $|I|\geq 1$ and
$\ann_{\ba}(p)=(0)$ for some $p\in\p$.  Let $q_{\p}$ be an
\ba-quadratic form on \ps and set $q=q_I\oplus q_{\p}$.  Let
$L=\eosp(q_{\infty})$ and $\ce=\eosp(q)$ (e.g. $L$ is a centreless
$D_J$- or $B_J$-graded Lie superalgebra for
$J=I\dot{\cup}\{\infty\}$, \cite{GN}).  Then
\[
    \der_{\bk}L=\ider L+\{d_{S,0}\in\de_{\ce}\mid
S\in\es_{\ce}\}
\]
and
\[
    \ider L\cap\{d_{S,0}\in\de_{\ce}\mid S\in\es_{\ce}\}=\ad\ce.
\]
Hence if $\es=\ce\rtimes\es_0$ then
\[
    \der_{\bk}L=\ider L\rtimes\{d_{S,0}\in\de_{\ce}\mid
S\in\es_0\}.
\]
In particular, if $\p=\n\oplus\calr$ with the \ba-quadratic form
$q=q_{\n}\oplus q_{\calr}$ where $\n$ and $\calr$ are free
\ba-supermodules with \ba-quadratic forms $q_{\n}$ and $q_{\calr}$
respectively such that $q_{\n}$ is almost diagonalizable and
$q_{\calr}\equiv 0$ then
\[
    \der_{\bk}L=\adloc L\rtimes\{d_{S,0}\in\de_{\ce}\mid
S\in\ed_A\calr\oplus(\der_KA)_{\m})
\]
where $d_{(\der_KA)_{\m},0}$ is isomorphic (as a Lie superalgebra)
to $\der_KA$.  In particular, if $L$ is of finite rank,
\[
    \der_{\bk}L=\ad L\rtimes\{d_{S,0}\in\de_{\ce}\mid
S\in\ed_A\calr\oplus(\der_KA)_{\m}\}.
\]
\end{co}
\textsl{Proof:}  Since $|I|\geq 1$, $H(I;\ba)$ contains a
two-dimensional invertibly diagonalizable free \ba-supermodule as
a direct summand and so the assertion holds by Proposition \ref{T
for even and q=inv diag+0}.  The rest follows
from Theorem \ref{der for q=inv diag+0}.\qed\\
\\
We can apply Theorem \ref{der for q=inv diag+0} to extended affine
Lie algebras of type $B$ and $D$:
\begin{co}
Let ${\cal{K}}$ be the centreless core of an extended affine Lie
algebra of type $B_l$ or $D_l$ ($l\geq 3,\ l\geq 4$ respectively)
with nullity $\nu$.  Then
\[
    \der_{\C}{\cal{K}}\cong
    \ider {\cal{K}}\rtimes\der_{\C}\C[t_1^{\pm
1},\ldots,t_{\nu}^{\pm 1}].
\]
\end{co}
The last assertion of the following was proved in
\cite[Proposition 2.3.4]{Kac2} and \cite[Proposition
3.1.2.3]{Scheunert} and is a special case of Theorem \ref{der for
q=inv diag+0}.
\begin{co}\label{kac and scheunert}
Let \bks be an algebraically closed field of characteristic 0.
\begin{itemize}
\item[(1)]  Let $m,n\in\Z^+$,
$q=q_{\{1,\ldots,m\}}\oplus\tilde{q}_{\{1,\ldots,n\}}$ over \bks
and let $q_0:\bk m_0\times\bk m_0\rightarrow \bk$ be the
\bk-quadratic form on a free supermodule $\bk m_0$ given by
$q(m_0,m_0)=1$.  Then $\eosp(q)=\osp(2m,2n)$ and $\eosp(q\oplus
q_0)=\osp(2m+1,2n)$. \item[(2)]  Let \bas be a superextension of
\bks and let $L=\osp(m,2n),\ m,n\in\Z^+$.  Then
$Der_{\bk}(L\otimes\ba)\cong \ider L\rtimes \der_{\bk}\ba$.  In
particular, for $\ba=\bk$, $L$ has no outer derivations.
\end{itemize}
\end{co}

\newpage


\end{document}